\documentclass[a4paper,10pt]{article}
\usepackage{amsfonts,amssymb,amsthm,amsmath}
\usepackage[matrix,arrow]{xy}

\theoremstyle{plain}
\newtheorem{theorem}{Theorem}[section]
\newtheorem*{Th*}{Theorem}
\newtheorem{lemma}[theorem]{Lemma}

\newtheorem*{Cor*}{Corollary}
\newtheorem{proposition}[theorem]{Proposition}

\theoremstyle{definition}

\theoremstyle{remark}

\numberwithin{equation}{section}

\usepackage{graphicx}
\usepackage{amsmath,amssymb}

\numberwithin{equation}{section}

\makeatletter
\def\Set@Scallop[#1]#2#3{{#1}\Parens{#2}{#3}}
\newcommand\DeclareScalableOperator[2]{%
  \expandafter\def\csname#1\endcsname{\@ifnextchar[{{#2}\Set@Scallop}{{#2}\Set@Scallop[{}]}}
}
\makeatother

\DeclareScalableOperator{Ct}{\mathcal{C}} 
\DeclareScalableOperator{Cc}{\mathcal{K}} 
\DeclareScalableOperator{Lp}{\mathbf{L}} 
\DeclareScalableOperator{Hom}{\mathrm{Hom}} 
\DeclareScalableOperator{Aut}{\mathrm{Aut}} 
\DeclareScalableOperator{Int}{\mathrm{Int}} 
\DeclareScalableOperator{aut}{\mathrm{aut}} 
\DeclareScalableOperator{End}{\mathrm{End}} 
\DeclareScalableOperator{GL}{\mathrm{GL}}
\DeclareScalableOperator{SLL}{\mathrm{SL}}
\DeclareScalableOperator{SO}{\mathrm{SO}}
\DeclareScalableOperator{sll}{\ger{sl}}
\DeclareScalableOperator{gl}{\ger{gl}}
\DeclareScalableOperator{co}{\mathrm{co}}

\newcommand\jBox[2]{{#1}\,\Box\,{#2}^*}

\DeclareMathOperator{\Ad}{\mathrm{Ad}}
\DeclareMathOperator{\ad}{\mathrm{ad}}
\DeclareMathOperator{\id}{\mathrm{id}}
\DeclareMathOperator\rk{\mathrm{rk}}
\DeclareMathOperator\tr{\mathrm{tr}}

\DeclareMathOperator{\im}{\mathrm{Im}}

\newcommand\var{\llcorner\kern-.3em\lrcorner}
\newcommand\clos{\overline}

\newcommand\ger{\mathfrak}

\newcommand{\fa}{for all }
\newcommand{\Fs}{For some }
\newcommand{\fs}{for some }
\newcommand\mathfa[1][{}]{\quad\text{\fa{#1} }}
\newcommand\smathfa[1][{}]{\ \text{\fa{#1} }}

\newcommand{\scth}{such that }
\newcommand{\AND}{and}

\newcommand\mathtxt[1]{\quad\text{{#1}}\quad}

\newcommand{\nd}{\mathtxt\AND}
\newcommand{\nda}{\ \text{\AND}\ }

\newcommand\vtheta{\vartheta}

\newcommand\reals{\mathbb{R}}
\newcommand\cplxs{\mathbb{C}}

\newcommand\vvoid{\varnothing}
\newcommand\sle{\leqslant}
\newcommand\sge{\geqslant}
\newcommand\limk{\lim\nolimits}

\makeatletter
\newcommand\Size[7][1]{
                                 \ifx#20%
                                        \def\r@l{}\def\r@m{}\def\r@r{}%
                                 \else%
                                    \ifx#21%
                                           \def\r@l{\bigl}\def\r@r{\bigr}\def\r@m{\bigm}%
                                    \else%
                                           \ifx#22%
                                                 \def\r@l{\Bigl}\def\r@r{\Bigr}\def\r@m{\Bigm}%
                                            \else%
                                                 \ifx#23%
                                                        \def\r@l{\biggl}\def\r@r{\biggr}\def\r@m{\biggm}%
                                                  \else
                                                        \ifx#24%
                                                              \def\r@l{\Biggl}\def\r@r{\Biggr}\def\r@m{\Biggm}%
                                                        \fi%
                                                  \fi%
                                            \fi%
                                      \fi%
                                 \fi%
                                 \ifx#10%
                                       \def\r@m{}%
                                 \fi%
                                 \r@l#3{#4}\r@m#5{#6}\r@r#7%
}%
\makeatother

\newcommand\Set[3]{
                                 \Size{#1}{\{}{#2}{|}{#3}{\}}%
}%
\newcommand\Scp[3]{
                                 \Size{#1}{(}{#2}{|}{#3}{)}%
}%
\newcommand\Rscp[3]{
                                 \Size[0]{#1}{(}{#2}{:}{#3}{)}%
}%
\newcommand\Parens[2]{
  \Size[0]{#1}{(}{#2}{}{}{)}
}
\newcommand\Bracks[2]{
  \Size[0]{#1}{[}{#2}{}{}{]}
}
\newcommand\Braces[2]{
  \Size[0]{#1}{\{}{#2}{}{}{\}}
}

\newcommand\Abs[2]{
  \Size[0]{#1}{\lvert}{#2}{}{}{\rvert}
}
\newcommand\Span[2]{
  \Size[0]{#1}{\langle}{#2}{}{}{\rangle}
}
\newcommand\trip[4]{
                                 \Size[0]{#1}{\{}{#2}{{{#3}^*}}{#4}{\}}%
}%

\makeatletter

\newif\if@smallmat
\newif\if@none
\newif\if@paren
\newif\if@brack
\newif\if@brace
\newif\if@vline
\newenvironment{Matrix}[2][1]
                                 {\ifx#20%
                                        \@smallmattrue%
                                  \else%
                                         \@smallmatfalse
                                  \fi%
                                  \ifx#11%
                                         \@nonefalse\@parentrue\@brackfalse\@bracefalse\@vlinefalse%
                                  \else%
                                       \ifx#12%
                                            \@nonefalse\@parenfalse\@bracktrue\@bracefalse\@vlinefalse%
                                        \else%
                                            \ifx#13%
                                                 \@nonefalse\@parenfalse\@brackfalse\@bracetrue\@vlinefalse%
                                            \else%
                                                 \ifx#14%
                                                       \@nonefalse\@parenfalse\@brackfalse\@bracefalse\@vlinetrue
                                                 \else%
                                                       \ifx#15%
                                                             \@nonefalse\@parenfalse\@brackfalse\@bracefalse\@vlinefalse%
                                                       \else%
                                                             \@nonetrue\@parenfalse\@brackfalse\@bracefalse\@vlinefalse%
                                                       \fi%
                                                 \fi%
                                            \fi%
                                        \fi%
                                   \fi%
                                   \if@smallmat%
                                        \if@none%
                                             \begin{smallmatrix}%
                                        \else%
                                            \if@paren%
                                                  \bigl(\begin{smallmatrix}%
                                            \else%
                                                  \if@brack%
                                                          \bigl[\begin{smallmatrix}%
                                                  \else%
                                                          \if@brace%
                                                               \bigl\{\begin{smallmatrix}%
                                                          \else%
                                                               \if@vline%
                                                                    \bigl\lvert\begin{smallmatrix}%
                                                                \else%
                                                                    \bigl\lVert\begin{smallmatrix}%
                                                                \fi%
                                                          \fi%
                                                  \fi%
                                            \fi%
                                        \fi%
                                   \else%
                                        \if@none%
                                             \begin{matrix}%
                                        \else%
                                            \if@paren%
                                                  \begin{pmatrix}%
                                            \else%
                                                  \if@brack%
                                                          \begin{bmatrix}%
                                                  \else%
                                                          \if@brace%
                                                               \begin{Bmatrix}%
                                                          \else%
                                                               \if@vline%
                                                                    \begin{vmatrix}%
                                                                \else%
                                                                    \begin{Vmatrix}%
                                                                \fi%
                                                          \fi%
                                                  \fi%
                                            \fi%
                                        \fi%
                                   \fi}%
                                  {\if@smallmat%
                                        \if@none%
                                             \end{smallmatrix}%
                                        \else%
                                            \if@paren%
                                                  \end{smallmatrix}\bigr)%
                                            \else%
                                                  \if@brack%
                                                          \end{smallmatrix}\bigr]%
                                                  \else%
                                                          \if@brace%
                                                               \end{smallmatrix}\bigr\}%
                                                          \else%
                                                               \if@vline%
                                                                    \end{smallmatrix}\bigr\rvert%
                                                                \else%
                                                                    \end{smallmatrix}\bigr\rVert%
                                                                \fi%
                                                          \fi%
                                                  \fi%
                                            \fi%
                                         \fi%
                                   \else%
                                        \if@none%
                                             \end{matrix}%
                                        \else%
                                            \if@paren%
                                                  \end{pmatrix}%
                                            \else%
                                                  \if@brack%
                                                          \end{bmatrix}%
                                                  \else%
                                                          \if@brace%
                                                               \end{Bmatrix}%
                                                          \else%
                                                               \if@vline%
                                                                    \end{vmatrix}%
                                                                \else%
                                                                    \end{Vmatrix}%
                                                                \fi%
                                                          \fi%
                                                  \fi%
                                            \fi%
                                        \fi%
                                   \fi}%

\makeatletter

\newcommand{\IfUpperCase}[1]{\begingroup 
  \protected@edef\@tempa{\expandafter\@firstofone\@firstofone#1.}%
  \expandafter\IfUpperCasE \@tempa\delimiter}

\def\IfUpperCasE #1#2\delimiter{%
  \protected@edef\@tempa{\meaning#1\meaning a}%
  \ifnum \expandafter\IfUppercaSE\@tempa \IfUppercaSE
   \endgroup \expandafter\@firstoftwo
  \else
   \endgroup \expandafter\@secondoftwo
  \fi}

\@ifundefined{strip@prefix}{\def\strip@prefix#1>{}}{}
\def\@tempa{the letter }
\edef\@tempa{\expandafter\strip@prefix\meaning\@tempa}

\expandafter\def\expandafter\IfUppercaSE\expandafter#\expandafter1\@tempa#2#3\IfUppercaSE{\uccode`#2=`#2 }

\newif\ifuc@se
\newcommand{\setuc@se}[1]{\IfUpperCase{#1}{\uc@setrue}{\uc@sefalse}}

\newcommand{\theoremn@me}[1]{\ifuc@se \lowercase{\csname#1name\endcsname}\ignorespaces%
  \else \edef\@temp{\lowercase{\lowercase{\csname#1name\endcsname}}}\@temp\ignorespaces%
  \fi}
\newcommand{\theoremn@mes}[1]{\ifuc@se \lowercase{\csname#1names\endcsname}\ignorespaces%
  \else \edef\@temp{\lowercase{\lowercase{\csname#1names\endcsname}}}\@temp\ignorespaces%
  \fi}
  
\newcommand\thmref[2]{\setuc@se{#1}\lowercase{{\theoremn@me{#1}\lowercase{\ref{#1:#2}}}}}
\newcommand\thmpartref[3]{\setcounter{the@rem}{#3}\setuc@se{#1}\lowercase{\Hyper@Link{#1:#2:\thethe@rem}{\theoremn@me{#1}\lowercase{\ref{#1:#2} \ref{#1:#2:\thethe@rem}}}}}
\newcommand\thmsref[3]{\setuc@se{#1}\lowercase{\theoremn@mes{#1}}\lowercase{\ref{#1:#2}}\ \AND\ \lowercase{\ref{#1:#3}}}
\newcommand\thmsthreeref[4]{\setuc@se{#1}\lowercase{\theoremn@mes{#1}}\lowercase{\ref{#1:#2}},\ \lowercase{\ref{#1:#3}}\ \AND\ \lowercase{\ref{#1:#4}}}
\newcommand\thmsfourref[5]{\setuc@se{#1}\lowercase{\theoremn@mes{#1}}\lowercase{\ref{#1:#2}},\ \lowercase{\ref{#1:#3}},\ \lowercase{\ref{#1:#4}}\ \AND\ \lowercase{\ref{#1:#5}}}
\def\thmsfiveref#1#2#3#4#5#6{\setuc@se{#1}\lowercase{\theoremn@mes{#1}}\lowercase{\ref{#1:#2}},\ \lowercase{\ref{#1:#3}},\ \lowercase{\ref{#1:#4}},\ \lowercase{\ref{#1:#5}}\ \AND\ \lowercase{\ref{#1:#6}}}

\newcommand{\DefTheorem}[2]{\newenvironment{#1}[1][\empty]{\ignorespaces\begin{#2}\ifx##1\empty{}\else\lowercase{\label{#1:##1}}\fi\ignorespaces}{\end{#2}\ignorespacesafterend}}

\DefTheorem{Th}{theorem}
\DefTheorem{Prop}{proposition}
\DefTheorem{Cor}{corollary}
\DefTheorem{Lem}{lemma}
\DefTheorem{Def}{definition}
\DefTheorem{Rem}{remark}
                         
\makeatother

\begin{document}

\title{Boundary Orbit Strata and Faces of Invariant Cones and Complex Ol'shanski\u\i{} Semigroups}

\author{Alexander Alldridge\thanks{A.~Alldridge was partially supported by the IRTG ``Geometry and Analysis of Symmetries'', funded by Deutsche Forschungsgemeinschaft (DFG), Minist\`ere de l'\'Education Nationale (MENESR), and Deutsch-Franz\"osische Hochschule (DFH-UFA).}}

\maketitle

\begin{abstract}
	Let $D=G/K$ be an irreducible Hermitian symmetric domain. Then $G$ is contained in a complexification $G_\cplxs$, and there exists a closed complex subsemigroup $G\subset\Gamma\subset G_\cplxs$ characterised by fact that all holomorphic discrete series representations of $G$ extend holomorphically to $\Gamma^\circ$, the so-called \emph{minimal Ol'shanski\u\i{}} semigroup.
	
	Parallel to the classical theory of boundary strata for the symmetric domain $D$, due to Kor\'anyi and Wolf, we give a detailed and complete description of the $K$-orbit type strata of $\Gamma$ as $K$-equivariant fibre bundles. They are given by the conjugacy classes of faces of the minimal invariant cone in the Lie algebra $\ger g$. 


\end{abstract}

\setcounter{section}{-1}
\section{Introduction}

\noindent
The boundary structure of Hermitian symmetric domains $D=G/K$ is well understood through the work of Pjatecki\u\i{}-Shapiro (for the classical domains), and of Wolf and Kor\'an\-yi, in their seminal papers from 1965 \cite{koranyi-wolf,wolf-koranyi}: Each of the strata is a $K$-equi\-va\-riant fibre bundle whose fibres are Hermitian symmetric domains of lower rank. This detailed understanding of the geometry of $\clos D$ has been fruitful, and is the basis of a variety of developments in representation theory, harmonic analysis, complex and differential geometry, Lie theory, and operator algebras. We mention a few developments. 

	The original motivation of Wolf--Kor\'anyi was to provide Siegel domain realisations for Hermitian symmetric domains, without recourse to their classification. The existence of such realisations alone has led to an extensive literature way beyond the scope of this introduction. 
		
	The study of compactifications of (locally and globally) symmetric spaces is of current and continued interest (we mention the recent monograph \cite{borel-ji}). As a prominent example, the Baily--Borel compactification of Hermitian symmetric domains has been studied intensively, with applications to moduli spaces of K3 surfaces, variation of Hodge structure, and modular forms, among others. Its understanding relies essentially on the Wolf--Kor\'anyi result. 
	
	The Wolf--Kor\'anyi theory has been generalised to complex flag manifolds \cite{wolf-orbit1,wolf-zierau} and thus played an important role in the realisation theory of tempered representation of semi-simple Lie group (compare the references in \cite{wolf-flagoverview}); recently, it has found applications to cycle spaces \cite{hong-huckleberry,wolf-zierau-cycle,wolf-cycle} and orbit duality in flag manifolds \cite{bremigan-lorch,matsuki-orbit3}. 
	
	Further applications of the original Wolf--Kor\'anyi theory include unitary highest weight representations \cite{inoue-hardy,davidson-stanke,arazy-upmeier}; Poisson integrals \cite{koranyi-poisson,johnson-koranyi,koufany-zhang,boussejra-koufany}; Hardy spaces on various domains \cite{chadli,olafsson-orsted-hardy,bertram-hilgert-geomhardy}; parahermitian or Cayley type symmetric spaces \cite{kaneyuki1,kaneyuki2}; Toeplitz operators \cite{upmeier,upmeier-index}. 
	
\medskip\noindent
In 1977, Gel'fand and Gindikin \cite{gel_gind} proposed to study the harmonic analysis of Lie groups of Hermitian type $G$ by considering them as extreme boundaries of certain complex domains in $G_\cplxs$, to which certain series of representations should extend holomorphically. This programme has been widely investigated; notably, it has led to the definition of the so-called Ol'shanski\u\i{} semigroups and to Hardy type spaces of holomorphic functions on their interiors \cite{olshanskii_holext,olshanskii_hardy,stanton}. More recent progress has been made through the study of so-called complex crowns \cite{kroetz-stanton1,kroetz-stanton2,kroetz-opdam}. 
	
	Although a great deal is known about Ol'shanski\u\i{} domains \cite{koufany-orsted,neeb-olshanskii,kroetz-hardy,kroetz-stein,kroetz-dual}, their boundary structure has as yet not been completely investigated. As in the case of Hermitian symmetric domains, we would expect that detailed and complete information on the $K$-orbit type strata ($K\subset G$ maximal compact) could lead to a better understanding of the geometry and analysis on these domains.
	
	To be more specific, fix an irreducible Hermitian symmetric domain $D=G/K$. The Lie algebra $\ger g$ of $G$ contains a minimal $G$-invariant closed convex cone $\Omega^-$ and $\Gamma=G\cdot\exp i\Omega^-\subset G_\cplxs$ is the minimal Ol'shanski\u\i{} semigroup. We describe all the faces of $\Omega^-$ (\thmref{Th}{facemincone}); each of them can be described explicitly, and gives rise to an Ol'shanski\u\i{} semigroup in the complexification of a certain subgroup of $G$. These subgroups are semidirect products $S\ltimes H$ where $S$ is the connected automorphism groups of a (convex) face of $D$ and $H$ is a certain generalised Heisenberg group related to the intersection of two maximal parabolic subalgebras of $\ger g$. 
	
	The relative interiors of the faces fall into $G$- (equivalently, $K$-) conjugacy classes, each of which forms exactly one of the $K$-orbit type strata of $\Omega^-$ (\thmref{Th}{faceconj}). One immediately deduces the $K$-orbit type stratification of $\Gamma$ (\thmref{Th}{olshanski-strata}). Each stratum is a $K$-equivariant fibre bundle whose fibres are the $G$-orbits of the `little Ol'shanski\u\i{}' semigroups alluded to above. In particular, the fibres are $K$-equivariantly homotopy equivalent to $K$ itself.
	
	While this result is in beautiful analogy to that of Wolf--Kor\'anyi, we stress that the structure of the Ol'shanski\u\i{} semigroups occuring as fibres is more complicated than that of $\Gamma$---their unit groups are not all Hermitian simple; rather, they have the structure of `generalised Jacobi groups'. Furthermore, we remark that all of the above statements can and will be made entirely explicit in the main text of this paper, by the use of the Jordan algebraic structure of the Harish-Chandra embedding of $D$. 
	
\medskip\noindent
Let us give a more detailed overview of our paper. In section 1, we collect several basic facts about symmetric domains, symmetric cones, and the associated Lie and Jordan algebraic objects. While most of the information we recall here can be easily extracted from the literature, some items are more specific. So, although this accounts for a rather lengthy glossary of results, we feel that it may serve as useful reference, in particular with regard to some of the more technical arguments of this article. 

Section 2 contains an account of the classification of nilpotent faces. In fact, in the course of the proof of the classification, we reprove the classification of conal nilpotent \emph{orbits}. Assuming the latter would not simplify our argument; indeed, our proof of the more precise result (\thmref{Th}{mincone-co-orbit}) is shorter than the existing proof of the classification of conal nilpotent orbits. The Theorem gives the description of all faces of the minimal (or maximal) invariant cone which contain a nilpotent element in their relative interior, and the decomposition of the nilpotent variety in the minimal cone into $K$-orbit type strata (which are the same as the conal nilpotent $G$-orbits). 

The main body of our work is the content of section 3. It culminates in the classification of the faces of the minimal invariant cone (\thmref{Th}{facemincone}), the characterisation of their conjugacy, and the description of the $K$-orbit type strata (\thmref{Th}{faceconj}). The basic observation is that each face generates a subalgebra (the \emph{face algebra}), and its structure is well understood due to the work of Hofmann, Hilgert, Neeb \emph{et al.} on invariant cones in Lie algebras. We show that the Levi complements of the face algebras are exactly the Lie algebras of complete holomorphic vector fields on the faces of the domain $D$. On the other hand, the centres and nilradicals of the face algebras can be understood through the classification of nilpotent faces. These considerations suffice to complete the classification of all those faces---both of the minimal and the maximal invariant cone---whose face algebra is \emph{non}-reductive. The classification in the case of \emph{reductive} face algebras only works well in the case of the minimal cone; it relies on the observation (\thmref{Lem}{extremalray-invcone}) that all extreme rays of the minimal cone are nilpotent (a fact which follows directly from the Jordan--Chevalley decomposition).

Finally, in section 4, we globalise the results of section 3 to the minimal Ol'shan\-sk\u\i{}i semigroup (\thmref{Th}{olshanski-strata}). Although the global results on the level of the semigroup are probably ultimately of greater interest than the infinitesimal results on the level of the minimal cone, the globalisation follows essentially by standard procedures. At this point, all the hard work has been done. 

\section{Bounded symmetric domains and\\ÊJordan triples}

	We begin with a revision of basic facts about bounded symmetric domains and related matters. We apologise to the reader for the tedium, but we will need the details. 

\subsection{Bounded symmetric domains, and their automorphism groups}

\paragraph{Bounded symmetric domains}
Let $Z$, $\dim Z=n$ be a complex vector space, $D\subset Z$ a circular bounded symmetric domain. Let $G$ be the connected component of $\Aut0D$. Then $G$ is a Lie group whose Lie algebra $\ger g$ of the set of \emph{complete} holomorphic vector fields on $D$. The bracket is $\Bracks1{h(z)\tfrac\partial{\partial z},k(z)\tfrac\partial{\partial z}}=\bigl(h'(z)k(z)-k'(z)h(z)\bigr)\tfrac\partial{\partial z}$ where $\Parens1{h(z)\tfrac\partial{\partial z}}f(z)=f'(z)h(z)$ for $h,f:D\to Z$. For $\xi\in\ger g$, $g\in\Aut0D$, $z\in D$, the adjoint action of $\Aut0D$ is $\Ad(g^{-1})(\xi)(z)=g'(z)^{-1}\xi(g(z))$. 

Let $\ger k$ be the set of \emph{linear} vector fields. The vector field $h_0=iz\tfrac\partial{\partial z}\in\ger z(\ger k)\setminus0$ generates the $\mathrm U(1)$-action. The isotropy $K$ of $G$ at $0$ hat the Lie algebra $\ger k$. Moreover, $D=G/K$; we have $Z(G)=1$; elements of $K$ are linear; $K$ is maximal compact and equals the fixed group of the Cartan involution $\vtheta=\Ad(-\id_D)$ \cite[Lemma~1.7]{upmeier_habil}, \cite[Sect.~1.2]{loos_boundsymmdom},\cite[Chapter~III, \S~7, Proposition~7.4]{helgason_dg_lg_sym}. 

\paragraph{Jordan triples} 
Consider the Cartan decomposition $\ger g=\ger k\oplus\ger p$; $\ger p\to Z:\xi\mapsto\xi(0)$ is a real linear isomorphism. Define $\xi_z^-\in\ger p$ by $\xi_z^-(0)=z$, further $Q_z(w)=z-\xi_z^-(w)$, $2Q_{u,w}=Q_{u+w}-Q_u-Q_w$; then $\{uv^*w\}=Q_{u,w}(v)$ is linear in $u$ and $w$, conjugate linear in $v$, and $\Parens0{\jBox uv}(w)=\trip0uvw$ satisfies \cite[Lemma~2.6]{loos_boundsymmdom}
\begin{equation}\label{eq:JTS}
\trip0uvw=\trip0wvu\ ,\ 
\Bracks1{\jBox uv,\jBox zw}=\jBox{\trip0 uvz}w-\jBox z{\trip0wuv}\ ,
\end{equation}
so $Z$ is a \emph{Jordan triple}. We have $\xi_u^-=(u-\trip0zuz)\tfrac\partial{\partial z}$ and identities \cite[Lemma 2.6.]{loos_boundsymmdom}
\begin{equation}\label{eq:fundid}
\Bracks1{\xi_u^-,\xi_v^-}=2\Parens0{\jBox uv-\jBox vu}\nd
\Bracks1{\Bracks1{\xi_u^-,\xi_v^-},\xi_w^-}=2\xi^-_{\trip0uvw-\trip0vuw}
\end{equation}
Here, we identify $\ger k$ with a subset of $\End0Z$. 

The trace form $\tr_Z\Parens0{\jBox uv}$ defines a positive Hermitian inner product on $Z$ \scth $\jBox vu=(\jBox uv)^*$, so the Jordan triple $Z$ is \emph{Hermitian} \cite[Lemma 2.6]{loos_boundsymmdom}, \cite{upmeier_bsymm}. W.r.t.~a certain norm, $D$ is the unit ball, and this sets up a bijection between isomorphism classes of Hermitian Jordan triples and of bounded symmetric domains \cite[Theorem~4.1]{loos_boundsymmdom}. Note that $D$ is convex, and $K$ is connected.

\paragraph{Triple automorphisms}ÊAny element $k\in\GL0Z$ such that $k(\jBox uv)k^{-1}=\jBox{ku}{(kv)}$ is called a \emph{triple automorphism}. $K$ is the connected component of the set $\Aut0Z$ of triple automorphisms \cite[Corollary~4.9]{loos_boundsymmdom}. The Lie algebra $\ger k$ of $\Aut0Z$ coincides with the set $\aut0Z$ of all \emph{triple derivations} ($\delta\in\End0Z$, $\Bracks1{\delta,\jBox uv}=\jBox{(\delta u)}v+\jBox u{(\delta v)}$). All triple derivations are inner, i.e.~$\ger k=\aut0Z=\Span0{\jBox uv-\jBox vu\mid u,v\in Z}_\reals$.

\subsection{Cartan decomposition and Killing form}

\paragraph{Polarisation of $\ger p$}

Denote the complexification of $\ger g$ by $\ger g_\cplxs$, etc. The decomposition $\ger g=\ger k\oplus\ger p$ gives $\ger g_\cplxs=\ger k_\cplxs\oplus\ger p_\cplxs$; $\vtheta$ extends to the conjugation of $\ger g_\cplxs$ w.r.t.~$\ger u=\ger k\oplus i\ger p$.

\begin{Lem}[pdecomp]
We have the vector space decomposition $\ger p_\cplxs=\ger p^+\oplus\ger p^-$ where
\[
\ger p^+=\Set1{u\tfrac\partial{\partial z}}{u\in Z}\nd\ger p^-=\Set1{\trip0zuz\tfrac\partial{\partial z}}{u\in Z}\ .
\]
Moreover, $\vtheta(\ger p^+)=\ger p^-$, and $\ger p^\pm$ are $\ger k_\cplxs$-invariant and Abelian.
\end{Lem}

\begin{proof}
We have 
\begin{equation}\label{eq:theta-pdecomp}
\vtheta\Parens1{u\tfrac\partial{\partial z}}=\tfrac12\vtheta(\xi_u^--i\xi_{iu}^-)=-\tfrac12(\xi_u^-+i\xi_{iu}^-)=\trip0zuz\tfrac\partial{\partial z}\ .
\end{equation}
The vector fields in $\ger p^+$ are constant, so $[\ger p^+,\ger p^+]=0$. Applying $\vtheta$ gives $[\ger p^-,\ger p^-\bigr]=0$. 

Any $\delta\in\ger k$ is linear, so $\Bracks1{\delta,u\tfrac\partial{\partial z}}=(\delta u)\tfrac\partial{\partial z}$ \fa $u\in Z$. Since $\ger k$ leaves $\ger u$ invariant and hence commutes with $\vtheta$, the assertion follows. 
\end{proof}

\begin{Lem}[kcentre-centraliser]
The centraliser of $h_0=iz\tfrac\partial{\partial z}$ in $\ger g$ is $\ger k$. More precisely, $\ad h_0=\pm i$ on $\ger p^\pm$.
\end{Lem}

\begin{proof}
Clearly, $\ad h_0=i$ on $\ger p^+$. The assertion follows from (\ref{eq:theta-pdecomp}). 
\end{proof}

\paragraph{Killing form} 

The Killing form $B$ of $\ger g$ is given by $B(\xi,\eta)=\tr_{\ger g}\Parens0{\ad\xi\ad\eta}$ \fa $\xi,\eta\in\ger g$. Its complex bilinear extension to $\ger g_\cplxs$ will also be denoted by $B$.

\begin{lemma}[{\cite[Lemma~4.2]{koecher_elementary}, \cite[Lemma~6.1]{upmeier_habil}}]\label{lem:killing}
The splitting $\ger g_\cplxs=\ger p^+\oplus\ger k_\cplxs\oplus\ger p^-$ is $B$-orthogonal, and $\ger p^\pm$ are $B$-isotropic. We have, \fa $\delta\in\aut0Z$, $u,v\in Z$, 
\begin{gather*}
	B(\delta,\jBox uv)=2\tr_Z\Parens1{\jBox{(\delta u)}v}\ ,\ B\Parens1{u\tfrac\partial{\partial z},\trip0zvz\tfrac\partial{\partial z}}=-4\tr_Z(\jBox uv)\ ,\\
	B\Parens1{\xi_u^-,\xi_v^-}=4\tr_Z(\jBox uv+\jBox vu)\ .
\end{gather*}
\end{lemma}

\begin{proof}
The decomposition is orthogonal since $\ger k_\cplxs$, $\ger p^\pm$ are distinct $\ad h_0$-eigenspaces. We have $\Bracks1{u\tfrac\partial{\partial z},\trip0zvz\tfrac\partial{\partial z}}=-2\cdot u\,\Box\,v^*$ from \eqref{eq:fundid}. Thus,
\[
B(\delta,\jBox uv)=-\tfrac12B\Parens1{\Bracks1{\delta,u\tfrac\partial{\partial z}},\trip0zvz\tfrac\partial{\partial z}}=-\tfrac12B\Parens1{(\delta u)\tfrac\partial{\partial z},\trip0zvz\tfrac\partial{\partial z}}\ .
\]
The 2nd equation implies the 1st; the 1st with $\delta=h_0$ implies the 2nd. But
\[
B\Parens0{h_0,\jBox uv}=i\tr_{\ger p^+}\ad(u\,\Box\,v^*)-i\tr_{\ger p^-}\ad(u\,\Box\,v^*)
\]
by \thmref{Lem}{kcentre-centraliser}. Moreover, by (\ref{eq:theta-pdecomp}) and $\vtheta(\jBox uv)=-\jBox vu$, 
\[
\Bracks1{\jBox uv,w\tfrac\partial{\partial z}}=\trip0uvw\tfrac\partial{\partial z}\ ,\ 
\Bracks1{\jBox uv,\trip0zwz\tfrac\partial{\partial z}}=-\trip1z{\trip0vuw}z\tfrac\partial{\partial z}\ .
\]
Because $\tr_Z\jBox uv=\overline{\tr_Z\jBox vu}$, we have $B\Parens0{h_0,\jBox uv}=2i\tr_Z(\jBox uv)$. 

The subspaces $\ger p^\pm$ are isotropic, since $B(\ger p^\pm,\ger p^\pm)=\mp iB\Parens0{h_0,[\ger p^\pm,\ger p^\pm]}=0$. Now, 
\[
B(\xi_u^-,\xi_v^-)=-B\Parens1{u\tfrac\partial{\partial z},\trip0zvz\tfrac\partial{\partial z}}-B\Parens1{\trip0zuz\tfrac\partial{\partial z},v\tfrac\partial{\partial z}}=4\tr_Z(\jBox uv+\jBox vu)\ .\]
\end{proof}

We remark that if $Z$ is a simple Jordan triple, then $\ger g$ is simple \cite[Theorem~4.4]{koecher_elementary}.

\subsection{Tripotents and faces of $D$}

\paragraph{Tripotents} 

An $e\in Z$ is a \emph{tripotent} if $\trip0eee=e$. The $Z_\lambda(e)=\ker(\jBox ee-\lambda)$ are called \emph{Peirce $\lambda$-spaces}. Then $Z=Z_0(e)\oplus Z_{1/2}(e)\oplus Z_1(e)$, orthogonal sum w.r.t.~the trace form \cite[Theorem~3.13]{loos_boundsymmdom}. We have the \emph{Peirce rules} \cite[Proposition~21.9]{upmeier_bsymm}
\[
\trip0{Z_\alpha(e)}{Z_\beta(e)}{Z_\gamma(e)}\subset Z_{\alpha-\beta+\gamma}(e)\ ,\ 
\trip0{Z_0(e)}{Z_1(e)}Z=\trip0{Z_1(e)}{Z_0(e)}Z=0\ .
\]
In particular, $Z_1(e)$ and $Z_0(e)$ are subtriples. For tripotents $e,c$, $\jBox ec=0$ if and only if $\trip0eec=0$ \cite[Lemma~3.9]{loos_boundsymmdom}; we call $e,c$ \emph{orthogonal} ($e\perp c$). Define an order by
\[
c\sle e\ :\Leftrightarrow\ \trip0{(e-c)}{(e-c)}{(e-c)}=e-c\nda c\perp e-c\ .
\]
Then non-zero minimal (maximal) tripotents are called \emph{primitive} (\emph{maximal}); $e$ is primitive (maximal) if and only if $Z_1(e)=\cplxs e$ ($Z_0(e)=0$). Any \emph{unitary} tripotent ($Z=Z_1(e)$) is maximal; the converse holds for $Z$ a Jordan algebra ($D$ of tube type). 

\paragraph{Frames, joint Pierce spaces}

A maximal orthogonal set $e_1,\dotsc,e_r$ of primitive tripotents is a \emph{frame}. In this case $r=\rk D$. Define $\rk Z=r$, $\rk e=\rk Z_1(e)$. Any tripotent equals $e_1+\dotsm+e_k$ for orthogonal primitive $e_j$ \cite[5.1, Theorem~3.11]{loos_boundsymmdom}. Given a frame, the \emph{joint Peirce spaces} are 
\[
Z_{ij}=\Set0{z\in Z}{\Braces0{e_k^{\vphantom*}e_k^*z}=\tfrac12(\delta_{ik}+\delta_{jk})\cdot z\ \forall k}\ ,\ 0\sle i\sle j\sle r\ .
\]
Then $Z=\bigoplus_{0\sle i\sle j\sle r}Z_{ij}$, $Z_{00}=0$, and $Z_{ii}=\cplxs e_i$ ($i>0$). If $Z$ is simple, $a=\dim Z_{ij}$, $b=\dim Z_{0j}$ are independent of $i,j$ and the frame, and $b=0$ exactly if $Z$ is a Jordan algebra. The \emph{canonical inner product} $\Scp0\var\var$ is the unique $K$-invariant inner product on $Z$ for which $\jBox vu=(\jBox uv)^*$ and $\Scp0 ee=1$ for every primitive tripotent $e$. Its restriction to any subtriple is canonical. For simple $Z$, $\Scp0 uv=\tfrac{2r}{2n-rb}\cdot\tr_Z(u\,\Box\,v^*)$.

\paragraph{Faces of $D$}

Given a convex set $C\subset\reals^d$, a subset $F\subset\clos C$ is a \emph{face} if any open line segment in $C$ intersecting $F$ lies in $F$. We let $F^\circ$ (the \emph{relative interior}) denote the interior of $F$ in its affine span. A hyperplane is \emph{supporting} if $C$ lies on one side of it. A proper face is \emph{exposed} if $F=\clos C\cap H$ for a supporting hyperplane $H$. 

For any tripotent $e\in Z$, define $D_0(e)=D\cap Z_0(e)$, the symmetric domain associated with $Z_0(e)$. Then $e+D_0(e)$ is a face of $D$, and this defines a bijection between tripotents of $Z$ and faces of $D$ \cite[Theorem~6.3]{loos_boundsymmdom}. All faces of $D$ are exposed. 

\begin{Def}[facial]
For any tripotent $e$, let $G_0(e)=\Aut[_0]0{D_0(e)}$. Then $G_0(e)=K_0(e)\cdot\exp\ger p_0(e)$ where $K_0(e)=\Aut[_0]0{Z_0(e)}$, $\ger p_0(e)=\Set0{\xi_u^-}{u\in Z_0(e)}$.
The Lie algebra of $K_0(e)$ is $\ger k_0(e)=\aut0{Z_0(e)}\subset\ger k$. In particular, $K_0(e)\subset K$, $G_0(e)$ is a closed subgroup of $G$, and $e+D_0(e)=G_0(e).e\cong G_0(e)/K_0(e)$. If $c\sle e$, then $G_0(e)\subset G_0(c)$. 
\end{Def}

\subsection{Symmetric cones and formally real Jordan algebras}

We conclude our preliminaries with a short section on symmetric cones. This may seem to be somewhat of a digression, but will be important in what follows. 

\paragraph{Unitary tripotents and Jordan algebras} 

A tripotent $e\in Z$ is \emph{unitary} if $Z=Z_1(e)$. It defines a composition $z\circ w=\trip0zew$ and an involution $z^*=\trip0eze$. Then $\circ$ is commutative, $e$ is a unit, and $z^2\circ(z\circ w)=z\circ(z^2\circ w)$, so $Z$ is a \emph{complex Jordan algebra} \cite[Proposition~13]{upmeier_bsymm}. The triple product is recovered by 
\begin{equation}\label{eq:tripvscircprod}
\{uv^*w\}=u\circ(v^*\circ w)-v^*\circ(w\circ u)+w\circ(u\circ v^*)\quad\text\fa\ u,v,w\in Z\ .
\end{equation}
The $\circ$-closed real form $X=\{x\in Z\mid x^*=x\}$ is a real \emph{Jordan algebra}. Furthermore, $x^2+y^2=0$ implies $x=y=0$, i.e.~$X$ is \emph{formally real} (or \emph{Euclidean}) \cite[Theorem~3.13]{loos_boundsymmdom}. Conversely, for formally real $X$, $X_\cplxs$ is a complex Jordan algebra whose underlying triple is Hermitian. The canonical real form of $Z_1(e)$ is denoted $X_1(e)$. 

\paragraph{Symmetric cones}

Let $X$, $\dim X=n$, be a real vector space with an inner product $\Rscp0\var\var$. A convex cone $\Omega\subset X$ is \emph{pointed} if it contains no affine line, \emph{solid} if its interior in $X$ is non-void, and \emph{regular} if it is both pointed and solid. The \emph{dual cone} by $\Omega^*=\Set1{x\in X}{\Rscp0x\Omega\sge 0}$ is pointed (solid) if and only if $\Omega$ is solid (pointed). Let $\Omega\subset X$ be a closed solid cone. Then $\GL0\Omega=\Set0{g\in\GL0X}{g\Omega=\Omega}$ is a closed subgroup of $\GL0X$; denote its Lie algebra by $\gl0\Omega$. $\Omega$ is \emph{symmetric} if $\Omega^*=\Omega$ and $\GL0\Omega$ acts transitively on $\Omega^\circ$. Any symmetric cone is pointed.

Assume $\Omega$ symmetric. Then $\vtheta(g)=(g^{-1})^t$ is a Cartan involution of the reductive group $\GL0\Omega$, with compact fixed group $\mathrm O(\Omega)=\mathrm O(X)\cap\GL0\Omega$, and we may fix $e\in\Omega^\circ$ such that the stabiliser $\GL0\Omega_e=\mathrm O(\Omega)$ \cite[Proposition~I.1.8]{faraut_koranyi}. Denote the Cartan decomposition $\gl0\Omega=\ger o(\Omega)\oplus\ger p(\Omega)$. The linear map $\xi\mapsto\xi(e):\ger p(\Omega)\to X$ is an isomorphism. Define $M_x\in\ger p(\Omega)$ by $M_x(e)=x$. Then $x\circ y=M_xy$ makes $X$ a formally real Jordan algebra with identity $e$ \cite[Theorem~III.3.1]{faraut_koranyi}. On the other hand, $\Omega=\{x^2\mid x\in X\}$. This sets up a bijection between isomorphism classes of symmetric cones and of formally real Jordan algebras \cite[Theorem III.2.1]{faraut_koranyi}. 

The connected component $\GL[_+]0\Omega$ of $\GL0\Omega$ is transitive on $\Omega^\circ$, and $\mathrm O(\Omega)$ is the set $\Aut0X$ of \emph{Jordan algebra automorphisms} ($k\in\GL0X$, $k(x\circ y)=(kx)\circ(ky)$) \cite[Theorem~III.5.1]{faraut_koranyi}. Its Lie algebra $\ger o(\Omega)$ is the set $\aut0X$ of all \emph{Jordan algebra derivations}  ($\delta\in\End0X$ \scth $\delta(x\circ y)=(\delta x)\circ y+x\circ(\delta y)$). It can be seen that $\Aut0X$ is the set of those triple automorphisms $k$ of $X\otimes\cplxs$ such that $ke=e$, and that $\aut0X$ consists of all triple derivations $\delta$ such that $\delta(e)=0$. 

\paragraph{Idempotents and Pierce decomposition}
Any $c\in X$ \scth $c^2=c$ is an \emph{idempotent}. Let $X_\lambda(c)=\ker(M_c-\lambda)$ is the \emph{Peirce $\lambda$-space}; we have the \emph{Peirce decomposition} $X=X_0(c)\oplus X_{1/2}(c)\oplus X_1(c)$, orthogonal w.r.t.~the trace form $\tr_X(M_{x\circ y})$. The trace form on $X$ is positive, symmetric, $\mathrm O(\Omega)$-invariant, and hence proportional to $\Rscp0\var\var$. 

As above, we define orthogonality and ordering of idempotents. The non-zero minimal idempotents are \emph{primitive}, and maximal orthogonal sets of primitive idempotents are \emph{frames}. Their common cardinality is $r=\rk X$. Then $\rk c=\rk X_1(c)=k$ if and only if $c=c_1+\dotsc+c_k$ for orthogonal primitive $c_j$. The \emph{canonical inner product} $\Scp0\var\var$ is the unique $\mathrm O(\Omega)$-invariant inner product for which $\Scp0{u\circ v}w=\Scp0 v{u\circ w}$ and $\Scp0cc=1$ for any primitive idempotent $c$. If $X$ is simple, then $\Scp0xy=\tfrac rn\cdot\tr_X\Parens0{M_{x\circ y}}$.

\paragraph{Orbits and faces of $\Omega$}

Let $X$ be a simple formally real Jordan algebra, and fix a frame $c_1,\dotsc,c_r$. Let $c^k=c_1+\dotsm+c_k$. The cone $\Omega$ of squares decomposes into the $r+1$ orbits $\GL[_+]0\Omega.c^k$, $k=0,\dotsc,r$ \cite[Proposition~IV.3.1]{faraut_koranyi}. With any idempotent $c$, we associate the cone of squares $\Omega_0(c)\subset X_0(c)=X_1(e-c)$. 

\begin{proposition}[{\cite{braun-koecher}}]\label{prop:symmetriccone-faces}
The set of faces of $\Omega$ consists of
\[
\Omega_0(c)=X_0(c)\cap\Omega=c^\perp\cap\Omega=\Set1{x^2}{x\in X_0(c)}\ ,\ c=c^2\in X\ .
\]
In particular, all the faces of $\Omega$ are exposed. The dual face of $\Omega_0(c)$ is $\Omega_0(e-c)$. Two faces $\Omega_0(c)$ and $\Omega_0(c')$ are $\GL[_+]0\Omega$-conjugate if and only if $\rk c=\rk c'$. 
\end{proposition}

\begin{proof}
The set of elements of rank $k$ in $\Omega$ is $\GL[_+]0\Omega.c^k$; hence the conjugacy \cite[Proposition~IV.3.1]{faraut_koranyi}. For $x\in X$, $x\in\Omega$ if and only if $M_x$ is positive semi-definite \cite[Proposition~III.2.2]{faraut_koranyi}, so $\Omega\cap X_0(c)=\Omega_0(c)$. Since $c\in\Omega=\Omega^*$, $c^\perp$ is a supporting hyperplane, and $c^\perp\cap\Omega$ is an exposed face. We have $c^\perp\supset\Omega_0(c)$. On the other hand, $\Omega\cap c^\perp\subset X_0(c)$ if $c^2=c$ \cite[Exercise~III.3]{faraut_koranyi}, so $\Omega_0(c)=c^\perp\cap\Omega$. 

More generally, $\Omega_0(e-c)=\Omega\cap\Omega_0(c)^\perp$, and $\Omega_0(e-c)$, $\Omega_0(c)$ are dual faces. The extreme rays of $\Omega$ are the $\Omega_0(c)$ where $\rk c=r-1$ \cite[Proposition~IV.3.2]{faraut_koranyi}. Since $\Omega$ is self-dual, any proper face $F\subsetneq\Omega$ has a non-trivial dual face. Hence, $\Omega_0(e-c)$ is a maximal proper face. The faces of $\Omega$ contained in $\Omega_0(e-c)$ are exactly the faces of $\Omega_0(e-c)$. The claim follows by induction.
\end{proof}

\section{Nilpotent orbits and faces, maximal parabolics, and principal faces}

We now return to our setting of a Hermitian Jordan triple $Z$ of dimension $n$ and the associated circular bounded symmetric domain $D\subset Z$. In this section, we introduce the minimal and maximal invariant cones in $\ger g$, and classify their nilpotent faces. On the way, we reprove the classification of conal nilpotent \emph{orbits}. We also introduce a class of faces (called \emph{principal}) which are associated to maximal parabolic subalgebras. 

\subsection{Weyl group-invariant cones}

Consider the positive symmetric invariant form defined by 
\begin{equation}\label{eq:invform}
\Rscp0\xi\eta=-B(\xi,\vtheta\eta)\mathfa\xi,\eta\in\ger g\ .
\end{equation}

\paragraph{Toral Cartan subalgebra} 

Fix a frame $e_1,\dotsc,e_r$ of $Z$. By \cite[Lemma~1.1-2]{upmeier_jordan_harm}, there exists a Cartan subalgebra $\ger t^{\vphantom+}=\ger t^+\oplus\ger t^-\subset\ger k$ where
\begin{equation}\label{eq:cptcsa}
	\ger t^-=\Span0{ie_j^{\vphantom*}\,\Box\,e_j^*\mid j=1,\dotsc,r}_\reals\nda \ger t^+=\Set0{\delta\in\ger t}{\delta e_j=0\ \smathfa j=1,\dotsc,r}\ .
\end{equation} 
By \thmref{Lem}{kcentre-centraliser}, $\ger t$ is a Cartan subalgebra of $\ger g$. Let $\ger t_\cplxs=\ger t\otimes\cplxs$, and $\Delta=\Delta(\ger g_\cplxs:\ger t_\cplxs)$ the associated root system. Let $h_0=iz\tfrac\partial{\partial z}$. By \thmref{Lem}{kcentre-centraliser}, for $\alpha\in\Delta$, 
\[
\ger g^\alpha_\cplxs\subset\ger k_\cplxs\ \Leftrightarrow\ \alpha\Parens0{h_0}=0 \nd\ger g^\alpha_\cplxs\subset\ger p_\cplxs\ \Leftrightarrow\ \alpha\Parens0{h_0}\neq0\ .
\]
This gives a partition of $\Delta$ into subsets $\Delta_c$ and $\Delta_n$ of 
\emph{compact} and \emph{non-compact} roots, respectively. We consider the Weyl groups $W=W(\Delta)$ and $W_c=W(\Delta_c)$. For $\alpha\in\Delta$, let $H_\alpha\in i\ger t$ be determined by the fact that $B(H_\alpha,\cdot)\in\reals\alpha$, and $\alpha(H_\alpha)=2$ \cite[ch.~IV]{knapp_repbook}.

\begin{Def}
	Let $\Phi$ be a positive system of $\Delta$. Let $\Phi_c=\Delta_c\cap\Phi$ and $\Phi_n=\Delta_n\cap\Phi$. The positive system $\Phi$ is \emph{adapted} \cite{neeb_holconv} if \fa $\alpha,\beta\in\Phi_n$, $\alpha+\beta\not\in\Delta$. Equivalently: Any $\Phi_c$-simple root is $\Phi$-simple; $\Phi_n$ is $W_c$-invariant; for some (any) total order on $\Span0\Delta_\reals$ defining $\Phi$, $\Phi_c<\Phi_n$; the set $\Delta_c\cup\Phi_n$ is parabolic \cite[Proposition~VII.2.12]{neeb_holconv}.
\end{Def}

\begin{Lem}[adapted-possys]
Let $\Delta_c^{++}\subset\Delta_c$ be a positive system, $\Delta_n^{++}=\Set0{\alpha\in\Delta}{-i\alpha\Parens0{h_0}>0}$ and $\Delta^{++}=\Delta_c^{++}\cup\Delta_n^{++}$. Then $\Delta^{++}$ is adapted, and $\ger p^\pm=\bigoplus_{\alpha\in\Delta_n^{++}}\ger g^{\pm\alpha}_\cplxs$.
\end{Lem}

\begin{Rem}
For simple $Z$, adapted positive systems are $\{\pm1\}\times W_c$-conjugate \cite[Lem\-ma~VII.2.16]{neeb_holconv}. Moreover, a non-compact simple Lie algebra has an adapted positive system if and only if it is Hermitian, if and only if it is the Lie algebra of complete holomorphic vector fields of a bounded symmetric domain \cite[Proposition~VII.2.14]{neeb_holconv}.
\end{Rem}

\paragraph{Minimal $W_c$-invariant cone} 
Consider the following polyhedral cones in $\ger t$,  
\[
\omega^-=\mathrm{cone}\,\langle iH_\alpha\mid\alpha\in\Delta^{++}_n\rangle\ ,\ 
\omega^+=\Set0{H\in\ger t}{-i\alpha(H)\sge0\ \smathfa\alpha\in\Delta_n^{++}}\ .
\]
Then $\omega^+=(-i\Delta^{++}_n)^*$ is the dual cone of $\omega^-$, and both cones are pointed and have non-empty interior. By \thmref{Lem}{adapted-possys}, $\omega^\pm$ are $W_c$-invariant. We have $\omega^-\subset\omega^+$ \cite[Lemma~10]{hc_semisimple_VI}. For $k=1,\dotsc,r$, define $\gamma_k\in i\ger t^*$ by $\gamma_k(e_\ell\,\Box\,e_\ell^*)=\delta_{k\ell}$ and $\gamma_k(\ger t^+)=0$. Then $(\gamma_k)$ is a strongly orthogonal set \cite[Lemma~1.3]{upmeier_jordan_harm}, i.e.~$\gamma_k\pm\gamma_\ell\not\in\Delta$, $k\neq\ell$. 

Since $iz\tfrac\partial{\partial z}\equiv\sum_{k=1}^r ie_k^{\phantom*}\Box\,e_k^*\ (\ger t^+)$, we find $\gamma_k\in\Delta_n^{++}$. There is a total vector space order on $i\ger t^*$ defining $\Delta^{++}$, such that $0<\gamma_1<\dotsm<\gamma_r$. Consequently, $\gamma_1,\dotsc,\gamma_r$ is the \emph{Harish-Chandra fundamental sequence} \cite[II.6]{hc_semisimple_VI}. In particular, $(\gamma_k)$ is a strongly orthogonal set of maximal cardinality \cite[Lemma~8 and Corollary]{hc_semisimple_VI}.

\begin{Def}
	A root $\alpha\in\Delta$ is \emph{long} if $\Abs0\alpha\sge\Abs0\beta$ \fa $\beta\in\Delta$ contained in the same irreducible subsystem of $\Delta$ as $\alpha$.\footnote{Any irreducible subsystem of $\Delta$ has at most two root lengths \cite[Chapter VI, \S~1.4, Proposition~12]{bourbaki}.}
\end{Def}

The $\gamma_k$ are long \cite[Theorem~2]{moore_compact2}, \cite[Lemma~1]{paneitz_detcones}. All positive, long non-compact roots lying in the same irreducible subsystem of $\Delta$ are $W_c$-conjugate \cite[Lemma~2]{paneitz_detcones}.

\begin{Lem}[extremalgens]
The extreme rays of $\omega^-$ are generated exactly by $iH_\alpha$, $\alpha\in\Delta_n^{++}$, $\alpha$ long. In particular, $\omega^-=\mathrm{cone}\,\langle i\sigma(e_j^{\phantom*}\Box\,e_j^*)\mid\sigma\in W_c\,,\,j=1,\dotsc,r\rangle$. 
\end{Lem}

\begin{proof}
By definition, the generators of the extreme rays of $\omega^-$ are among the $H_\alpha$, $\alpha\in\Delta_n^{++}$. Since $\Delta$, and hence $\omega^-$, decomposes according to the decomposition of $\ger g$ into simple factors, we may assume w.l.o.g.~that $\ger g$ be simple. For any short $\gamma\in\Delta_n^{++}$, $\gamma=\frac12(\gamma_k+\gamma_\ell)$, some $ k\neq\ell$ \cite[Lemma~1]{paneitz_detcones}. Hence, $4\Abs0\gamma^2=\Abs0{\gamma_k}^2+\Abs0{\gamma_\ell}^2=2\Abs0{\gamma_k}^2$, and
\[
\Rscp0{H_\gamma}\xi=2\Abs0\gamma^{-2}\gamma(\xi)=2\Abs0{\gamma_k}^{-2}\gamma_k(\xi)+2\Abs0{\gamma_\ell}^{-2}\gamma_\ell(\xi)=\Rscp0{H_{\gamma_k}+H_{\gamma_\ell}}\xi
\]
\fa $\xi\in\ger t$. Hence, $iH_\gamma=iH_{\gamma_k}+iH_{\gamma_\ell}$ lies in the interior of a face of dimension at least $2$. On the other hand, $\omega^-$ being polyhedral, there is $\alpha\in\Delta_n^{++}$, necessarily long, such that $i\reals_{\sge0}\cdot H_\alpha$ is extreme; but all such $iH_\alpha$ are $W_c$-conjugate.

Returning to the semi-simple case, by maximality, any irreducible factor of $\Delta$ contains some $\gamma_k$. Moreover, any long $\alpha\in\Delta_n^{++}$ is $W_c$-conjugate to any $\gamma_k$ contained in the same irreducible factor. By \thmref{Lem}{killing} and (\ref{eq:cptcsa}), $e_j\,\Box\,e_j^*$ is proportional to $H_{\gamma_j}$. Suffices now to note that any polyhedral cone is generated by its extreme rays.
\end{proof}

\begin{Lem}[longroots]
Let $\gamma\in\Delta_n^{++}$ be long. There is a frame $c_1,\dotsc,c_r$ such that $\ger t$ is given by (\ref{eq:cptcsa}) (for $e_j=c_j$), and an integer $\ell$ \scth $\gamma(c_k^{\phantom*}\Box\,c_k^*)=\delta_{k\ell}$ and $\gamma(\ger t^+)=0$.
\end{Lem}

\begin{proof}
\Fs $\ell$, $\gamma_\ell$ and $\gamma$ lie in the same irreducible factor of $\Delta$; there is some $\sigma\in W_c$ such that $\sigma\gamma_\ell=\gamma$. Then $\sigma=\Ad(k)$ for some $k\in N_K(\ger t)$ \cite[Theorem~4.54]{knapp_beyond}. Since $k\in\Aut0Z$, the $ke_j$, $j=1,\dotsc,r$, are orthogonal primitive tripotents, and 
\[
\Ad(k)\Parens1{e_j^{\phantom*}\Box\,e_j^*\tfrac\partial{\partial z}}=k^{-1\prime}(z)^{-1}\{e_j^{\phantom*}e_j^*k^{-1}(z)\}\tfrac\partial{\partial z}=k(e_j^{\phantom*}\Box\,e_j^*)k^{-1}=(ke_j^{\phantom*})\,\Box\,(ke_j)^*
\]
where we identify linear maps and linear vector fields. Since $\Ad(k)$ normalises $\ger t$, we have a decomposition as stated. By the definition of $\gamma_\ell$, the lemma follows.
\end{proof}

\begin{Cor}[extremalgens-tripotents]
The extreme rays of $\omega^-$ are generated by the $i\cdot\jBox ee$ where $e$ is a primitive tripotent $W_c$-conjugate to an element of the frame $e_1,\dotsc,e_r$. 
\end{Cor}

\paragraph{Relation to the Weyl chamber}

The Weyl chamber associated to $\Delta^{++}$ is 
\[
c_+=\Set1{H\in\ger t}{-i\alpha(H)>0\ \smathfa\alpha\in\Delta^{++}}\ .
\]
By definition, it is obvious that $c_+\subset\omega^{+\circ}$. In fact, $\clos{c_+}$ is a fundamental domain for the action of $W_c$ on $\omega^+$ \cite[Lemma~I.5]{neeb-convcoadj}. From this, one immediate deduces the following statement.

\begin{Lem}
	Let $\Pi=(\alpha_k)$ be the simple system defining $\Delta^{++}$, and define $\omega_k\in\ger t$ by $\alpha_k(\omega_\ell)=i\delta_{k\ell}$. Then the generators of extreme rays of $\omega^+$ belong to $\bigcup_kW_c.\omega_k$. 
\end{Lem}

\subsection{Minimal and maximal invariant cones}

From now on, we assume that $Z$ be simple. Then $\ger z(\ger k)=\reals\cdot h_0$ where $h_0=iz\frac\partial{\partial z}$. 

\paragraph{Maximal and minimal cone} Consider the map $\Omega\mapsto\omega=\Omega\cap\ger t$ from the set of closed pointed $G$-invariant convex cones $\Omega\subset\ger g$ with non-trivial interior to the set of closed $W_c$-invariant convex cones $\omega$ \scth $\omega^-\subset\omega\subset\omega^+$. It is an order-preserving bijection \cite[Theorem~2]{paneitz_detcones}, and $\Omega=\{\xi\in\ger g\mid p_{\ger t}(\mathcal O_\xi)\subset\omega\}$ where $\mathcal O_\xi=\Ad(G)(\xi)$ and $p_{\ger t}$ is the orthogonal projection onto $\ger t$. Moreover, $\Omega^*\cap\ger t=(\Omega\cap\ger t)^*$ \cite[Theorem~3]{paneitz_detcones} and any orbit in $\Omega^\circ$ intersects the relative interior of $\omega$ non-trivially. 

Let $\Omega^-$ be the closed $G$-invariant convex cone generated by $iz\frac\partial{\partial z}$. Then we have $iz\,\frac\partial{\partial z}\in\omega^-\cap\Omega^{-\circ}$ \cite[Lemma~3]{paneitz_detcones}. All invariant cones with non-void interior have a $K$-fixed vector \cite[\S~2]{vinberg}, so $\Omega^-$ is minimal among invariant cones with non-void interior, and its dual $\Omega^+$ is maximal among pointed invariant cones.\footnote{For the invariance of $\Omega^+$, observe that $\Rscp0{\Ad(g)(x)}y=\Rscp0x{\Ad(\vtheta(g))(y)}$ \fa $x,y\in\ger g$, $g\in G$.} From this, it follows that $\Omega^\pm\cap\ger t=\omega^\pm$. The following result clarifies the structure of the set of semi-simple elements contained in $\Omega^\pm$.

\begin{Prop}[cone-regular]
Let $\xi\in\Omega^\pm$ be semi-simple. Then $\xi$ is conjugate to an element of $\omega^\pm$. If, in addition, $\xi$ is regular, then $\xi$ is conjugate to an element of $\omega^{+\circ}$ and hence contained in $\Omega^{+\circ}$. 
\end{Prop}

\begin{proof}
The orbit $\mathcal O_\xi$ is closed \cite[Proposition~1.3.5.5]{warner_semsimlg}. Hence, it intersects $\ger t$ \cite[Theorem~5.11]{hnp-sympconvcoadj}. This proves the first statement. The second now follows immediately from the fact that the set of regular semi-simple elements is open \cite[Proposition~1.3.4.1]{warner_semsimlg}, and that the centraliser of $\xi$ is a compact Cartan subalgebra. 
\end{proof}

\subsection{Tripotents, nilpotent faces, and nilpotent orbits of convex type}

Although Cayley triples have been extensively studied in the literature, we have to redo some of their theory to derive our result. In particular, we are interested in the following subclass of Cayley triples. 

\paragraph{$(H_1)$-Cayley triples}

A Lie algebra $\ger a$ is \emph{quasihermitian}, if $\ger b=\ger z_{\ger a}(\ger z(\ger b))$ for some maximal compact subalgebra $\ger b\subset\ger a$ containing a Cartan subalgebra of $\ger a$. If $\ger a$ is simple and non-compact, it is called \emph{Hermitian} if some maximal compact subalgebra has non-trivial centre. A reductive Lie algebra $\ger a$ is quasihermitian if and only if it is the direct sum of a maximal compact ideal and Hermitian simple ideals \cite{neeb_invconv}.

Consider the basis of $\sll0{2,\reals}$ given by $H=\begin{Matrix}0
1&0\\0&-1
\end{Matrix}$, $X^+=
\begin{Matrix}0
0&1\\0&0
\end{Matrix}$, $X^-=\begin{Matrix}0
0&0\\1&0
\end{Matrix}$. Let $\ger h$ be a quasihermitian reductive Lie algebra. Recall that $(h,x^+,x^-)\in\ger h^3$ is called an \emph{$\ger{sl}_2$-triple} if the associated linear map, defined by $H\mapsto h$ and $X^\pm\mapsto x^\pm$, is a Lie algebra monomorphism. $x^+$ is called the \emph{nilpositive} element of the triple. Given a Cartan involution $\theta$, an $\ger{sl}_2$-triple $(h,x^+,x^-)$ is called a \emph{Cayley triple} if $\theta(x^+)=-x^-$. 

An element $h_0\in\ger h$ is called an \emph{$H$-element} if $\ger z_{\ger h}(h_0)$ is maximally compactly embedded and $\mathop{\mathrm{Sp}}(\ad h_0)=\{0,\pm i\}$. Any $H$-element is semi-simple. With any $H$-element $h_0$, there is associated a unique Cartan involution $\theta=2\ad(h_0)^2+1$. For example, $iz\frac\partial{\partial z}$ is an $H$-element of $\ger g$, and $Z=\frac12(X^+-X^-)$ is an $H$-element of $\ger{sl}_2$. A homomorphism $\phi:\ger h\to\ger h'$ of quasihermitian reductive Lie algebras with fixed $H$-elements $h_0\in\ger h$ and $h_0'\in\ger h'$ is called an $(H_1)$-homomorphism if $\ad h_0'\circ\phi=\phi\circ\ad h_0$. 

Given an $H$-element $h_0\in\ger h$ with associated Cartan involution $\theta$, any Cayley triple $(h,x^+,x^-)$ will be called an \emph{$(H_1)$-Cayley triple} if the associated homomorphism $\ger{sl}_2\to\ger h$ is an $(H_1)$-homomorphism (relative to the $H$-elements $Z$ and $h_0$).

\begin{Lem}[cayley-triple]
	Let $h_0$ be an $H$-element in the quasihermitian reductive Lie algebra $\ger h$ with associated Cartan involution $\theta$, and let $x\in\ger h$. Then $x=x^+$ for some Cayley triple $(h,x^+,x^-)$ if and only if the following equation holds:
	\begin{equation}\label{eq:cayleytriple}
		[[\theta(x),x],x]=2x\ .
	\end{equation}
	This Cayley triple is unique. In this case, $\Span0{h,x^\pm}_\reals$ is $\ad h_0$-invariant if and only if $[h_0,x]=\pm\frac12[\theta(x),x]=\pm\frac12h$, and the triple is $(H_1)$ if and only if the sign is $+$. 
\end{Lem}

\begin{proof}
	If $x=x^+$ for some Cayley triple $(h,x^+,x^-)$, then $h=[x^+,x^-]=-[x,\theta(x)]$ and of course $x^-=-\theta(x)$. In particular, $(h,x^+,x^-)$ is unique, and (\ref{eq:cayleytriple}) holds.

	If equation (\ref{eq:cayleytriple}) holds, we \emph{define} $x^+=x$, $x^-=-\theta(x)$, $h=-[x,\theta(x)]$. Then we have $\theta(h)=-\theta([x,\theta(x)])=[x,\theta(x)]=-h$, 
	\[
		[h,x]=2x\nd [h,y]=\theta([h,x])=2\theta(x)=-2y\ .
	\]
	Thus, in this case, $(h,x^+,x^-)$ is a Cayley triple. 
	
	Next, observe that $\ger z_{\ger h}(h_0)=\ker(1-\theta)$. Setting $z=\tfrac12(x^+-x^-)$, this implies that $[h_0,z]=0$, so that $\ad(h_0)$ leaves the eigenspaces of $\ad z$ invariant. We have
	\[
		[z,h\pm i(x^++x^-)]=\pm i(h\pm i(x^++x^-))\ .
	\]
	If $\ad(h_0)$ leaves $\Span0{h,x^\pm}_\reals$ invariant, this implies that $h\pm i(x^++x^-)$ is an eigenvector of $\ad h_0$, for the eigenvalue $i$ or $-i$. The triple is $(H_1)$ if and only if the sign of the eigenvalue is the same as for $\ad z$. Moreover, again because $\ker(1-\theta)$ centralises $h_0$, $[h_0,x]=\frac12[h_0,x-\theta(x)]$, so $[h_0,x]=-[h_0,\theta(x)$, and $2[h_0,x]=[h_0,x^++x^-]$. Taking imaginary parts in the eigenvalue equation, $[h_0,x]=\pm\frac12[\theta(x),x]$, and the ($H_1$) condition amounts to the requirement that the sign be $+$.
\end{proof}

\begin{Prop}[cayleyh1]
	Fix the $H$-element $h_0=iz\frac\partial{\partial z}$, and let $x\in\ger g\setminus0$ be the nilpositive element of some Cayley triple. This triple is $(H_1)$ if and only if $x\in\Omega^-$, if and only if $x\in\Omega^+$. In particular, the nilpotent elements of $\Omega^+$ belong to $\Omega^-$. 
\end{Prop}

\begin{proof}
	Let $(h,x^+,x^-)$ be the Cayley triple with $x=x^+$. If $\Span0{h,x^\pm}_\reals$ is $\ad h_0$-invariant, then $[h_0,x]=\pm\frac12h$ by \thmref{Lem}{cayley-triple}. Thus we compute $e^{t\ad(x)}(h_0)=h_0\mp\tfrac t2h\pm\tfrac{t^2}2x$, and $\pm x=\lim_{t\to\infty}2t^{-2}e^{t\ad(x)}(h_0)\in\Omega^-$. 
	
	By \thmref{Lem}{cayley-triple}, if the triple is ($H_1$), then $x\in\Omega^-$. If $\Span0{h,x^\pm}_\reals$ is $\ad h_0$-invariant, let $x\in\Omega^+$, and assume that the triple is not $(H_1)$. Then $[h_0,x]=-\tfrac12[\vtheta(x),x]$ and $-x\in\Omega^-$. But $\Omega^-\subset\Omega^+$, and $\Omega^+$ is pointed. This is a contradiction, so the triple must be $(H_1)$. 
	
	We need to check that $x\in\Omega^+$ implies that $\Span0{h,x^\pm}_\reals$ is $\ad h_0$-invariant. It is sufficient to prove that $u^+=h+i(x^++x^-)\in\ger p^+\cup\ger p^-$. Up to $K$-conjugacy, we may assume that $x^+-x^-\in\ger t$. Since $-x^-=\vtheta(x)\in\Omega^+$, we have $z=\tfrac12(x^+-x^-)\in\omega^+$, so $-i\alpha(z)\sge0$ \fa $\alpha\in\Delta_n^{++}$. Since $[z,u^+]=iu^+$, we see that $z\in\ger p^+=\bigoplus_{\alpha\in\Delta_n^{++}}\ger g_\cplxs^\alpha$. Hence, $x\in\Omega^+$ implies that $\Span0{h,x^\pm}_\reals$ is $\ad h_0$-invariant. 
	
	Finally, any nilpotent element is $G$-conjugate to a nilpotent element belonging to a Cayley triple \cite[Theorems 9.2.1, 9.4.1]{collingwood-mcgovern}, so the claim follows. 
\end{proof}

For $u\in Z$, define the \emph{Cayley vector field} $\xi_u^+=-i\xi_{iu}^-=(u+\trip0zuz)\tfrac\partial{\partial z}$, and
\[
X_u^\pm=\tfrac12\Parens1{\xi_{-iu}^-\pm\tfrac12[\xi_u^-,\xi_{-iu}^-]}=\tfrac12(\xi_{-iu}^-\pm 2i\jBox uu)\ .
\]
For later use, we record the following simple formula:
\begin{equation}\label{eq:kaction-nilpotent}
	\Ad(k)(X_u^\pm)=\tfrac12\xi_{-iku}^-\pm\Ad(k)(i\jBox uu)=\tfrac12\xi_{-iku}^-\pm i\jBox{(ku)}{(ku)}=X_{ku}^\pm\ .
\end{equation}

\begin{Prop}[sl2triple]
Let $e,c\neq0$ be tripotents and $\ger s^e=\Span0{\xi_e^-,X_e^\pm}_\reals$. Then $(\xi_e^-,X_e^+,X_e^-)$ is an $(H_1)$-Cayley triple and $\pm X_e^\pm\in\Omega^-$. Moreover, $[\ger s^e,\ger s^c]=0$ if only if $e\perp c$.
\end{Prop}

\begin{proof}
First, note $[[\xi_{ae}^-,\xi_{be}^-],\xi_{ce}^-]=4\im(a\bar b)\xi_{ic\cdot e}^-$ \fa $a,b,c\in\cplxs$, whence
\[
[\xi_e^-,X_e^\pm]=\tfrac12[\xi_e^-,\xi_{-ie}^-]\mp\xi_{ie}^-=\pm2X_e^\pm\ ,\ 
[X_e^+,X_e^-]=\tfrac14[[\xi_e^-,\xi_{-ie}^-],\xi_{-ie}^-]=\xi_e^-\ .
\]
Clearly, $X_e^-=\vtheta(X_e^+)$, and $\Bracks1{iz\tfrac\partial{\partial z},X_e^-}=\frac12\Bracks1{iz\tfrac\partial{\partial z},\xi_{-ie}^-}=\frac12\xi_e^-$. Hence, the triple $(\xi_e^-,X_e^+,X_e^-)$ is an $(H_1)$-Cayley triple, and $\pm X_e^\pm\in\Omega^-$ by \thmref{Lem}{cayley-triple} and \thmref{Prop}{cayleyh1}. 

Next, observe $\xi_a^\pm\in\ger s^a$ for $a=e,c$. Since
\[
[\xi_e^-,\xi_c^-\bigr]-[\xi_e^-,\xi_c^+]=[\xi_e^-,\xi_c^-+i\xi_{ic}^-]=4\jBox ec\ ,\ [\xi_e^-,\xi_c^-]=[\jBox ee,\jBox cc]\ ,
\]
$e$ and $c$ are orthogonal if and only if $[\ger s^e,\ger s^c]=0$.
\end{proof}

\begin{Rem}
Paneitz \cite[Lemma~4]{paneitz_detcones} proves that $X_e^+\in\Omega^-$ for $e$ primitive.
\end{Rem}

\begin{Prop}[discembed-tripotent]
Let $(h,x^+,x^-)$ be an $(H_1)$-Cayley triple. Then there exists a unique non-zero tripotent $e\in Z$ such that $h=\xi_e^-$ and $x^\pm=X_e^\pm$.
\end{Prop}

\begin{proof}
We have $h\in\ger p$, so $h=\xi_e^-$ \fs $e\in Z\setminus0$. Set $z=\frac12(x^+-x^-)\in\ger k=\aut0Z$. The value $z(e)\in Z$ makes sense, and $\xi_{z(e)}^-=[z,h]=-(x^++x^-)$. 

By assumption, $\ad z$ and $\ad h_0$ (where $h_0=iz\frac\partial{\partial z}$) coincide on $\cplxs\Span0{h,x^+,x^-}_\reals$. Thus, $\xi_{z(e)}^-=[h_0,\xi_e^-]=\xi_{ie}^-$. This shows that $z(e)=ie$. Next,
\[
	\xi_{i\trip0eee}^-=\tfrac14[[\xi_e^-,\xi_{-ie}^*],\xi_e^-]=\tfrac14[[h,x^++x^-],h]=[z,h]=\xi_{ie}^-\ ,
\]
so $e=\trip0eee$. We have $x^++x^-=\xi_{ie}^-$, $x^+-x^-=\tfrac12[h,x^++x^-]=\tfrac12[\xi_e^-,\xi_{-ie}^-]$. 
\end{proof}

\begin{Rem}
	The result \cite[Proposition~4.1]{satake} seems to be somewhat similar.
\end{Rem}

We now introduce certain Heisenberg algebras associated to tripotents of $Z$. They will play a major role in the determination and description of the faces of $\Omega^\pm$.

\paragraph{Conal Heisenberg algebras}

In what follows, $e,c$ shall denote tripotents. 

\begin{Def}
Given $e$, and any set $A\subset\cplxs$, let $\ger g^e[A]=\bigoplus_{\lambda\in A}\ker(\ad\xi_e^--\lambda)$. Then $\ger g=\ger g^e[-2,-1,0,1,2]$ by \cite[Lemma~9.14]{loos_boundsymmdom}. Moreover (\emph{loc.cit.}), 
\[
\ger g^e[0]=\ger k^e\oplus\Set0{\xi_u^-}{u\in Z_0(e)\oplus X_1(e)}\mathtxt{where}\ger k^e=\Set0{\delta\in\ger k}{\delta e=0}\ .
\]
Let 
\begin{gather*}
\eta_u^e=\xi_u^-+[\xi_e^-,\xi_u^-]=\xi_u^-+2(\jBox eu-\jBox ue)\\
\zeta_u^e=\xi_u^-+\tfrac12[\xi_e^-,\xi_u^-]=\xi_u^-+(\jBox eu-\jBox ue)\mathfa u\in Z\ .
\end{gather*}
Then $X_e^\pm=\frac12\zeta_{\mp ie}^{\pm e}$, and (\emph{loc.cit.})
\[
\ger g^e[\pm 1]=\Set0{\eta_u^{\pm e}}{u\in Z_{1/2}(e)}\ ,\ \ger g^e[\pm 2]=\Set0{\zeta_u^{\pm e}}{u\in iX_1(e)}\ .
\]

Furthermore, $\ger q^e=\ger g^e[0,1,2]$ is a maximal parabolic \cite[Proposition~9.21]{loos_boundsymmdom}, and $\ger h^e=\ger g^e[1,2]$ is its nilradical. We call $\ger h^e$ a \emph{conal Heisenberg algebra}.

Recall that $\ger k_0(e)=\aut0{Z_0(e)}$ and $\ger k_1(e)=\aut0{X_1(e)}$ are, respectively, the set of triple derivations of $Z_0(e)$, and the set of algebra derivations of $X_1(e)$. Similarly, we consider $\ger p_0(e)=\Set0{\xi_u^-}{u\in Z_0(e)}$ and $\ger p_1(e)=\Set0{\xi_u^-}{u\in X_1(e)}$.

We already know that $\ger g_0(e)=\ger k_0(e)\oplus\ger p_0(e)$ is the set of complete holomorphic vector fields on $B_0(e)$. Let $\ger g_1(e)=\ger k_1(e)\oplus\ger p_1(e)$. Then by \cite[Lem\-ma~21.16]{upmeier_bsymm}, 
\[
\Ad(\gamma_e)(\xi_u^-)=2M_u\nd\Ad(\gamma_e)(\delta)=\delta\mathfa u\in X_1(e)\,,\,\delta\in\aut0{X_1(e)}
\]
where $M_u(v)=u\circ v$, so that $\Ad(\gamma_e)(\ger g_1(e))=\gl0{\Omega_1(e)}$.

We have $[\ger g_0(e),\ger g_1(e)]=0$ by the Peirce rules. Let $\ger m^e=\ger k^e\cap(\ger k_0(e)\oplus\ger k_1(e))^\perp$. If we let $\ger a=\Span0{\xi_{e_1}^-,\dotsc,\xi_{e_r}^-}_\reals$ for some frame such that $e_j\sle e$ or $e_j\perp e$ \fa $j$, then $\ger m^e\subset\ger z_{\ger k}(\ger a)$. Using this fact, it is easy to see that $\ger m^e$ leaves $\ger g_i(e)$ ($i=0,1$) invariant, so that, as Lie algebras,
\begin{equation}\label{eq:zerograding-decomp}
	\ger g^e[0]=\ger g_0(e)\oplus\ger g_1(e)\oplus\ger m^e=\aut0{D_0(e)}\oplus\Ad(\gamma_e^{-1})(\ger{gl}(\Omega_1(e))\oplus\ger m^e\ .
\end{equation}

Define a linear isomorphism $\phi^e:Z_{1/2}(e)\oplus X_1(e)\to\ger h^e$ by 
\begin{equation}\label{eq:heisenbergiso}
\phi^e(u,v)=\eta^e_u+\zeta_{-iv/2}^e\ .
\end{equation}
\end{Def}

\begin{Def}
	Let $U$, $V$ be complex vector spaces, $V$ be endowed with an antilinear involution ${}^*$, and $K$ be a closed convex cone \scth $x^*=x$ \fa $x\in K$. A sesquilinear map $\phi:U\times U\to V$ \scth $\phi(u,v)^*=\phi(v,u)$ and $\phi(u,u)\in K\setminus0$ \fa $u\neq0$ is called $K$-\emph{positive Hermitian}. 
\end{Def}

\begin{Prop}[heisenberg-alg]
	Define $h_e:Z_{1/2}(e)\times Z_{1/2}(e)\to Z_1(e)$ by $h_e(u,v)=8\cdot\trip0uve$, and $q_e(u,v)=\im h_e(u,v)=4i\cdot\Parens0{\trip0vue-\trip0uve}\in X_1(e)$. Then $h_e$ is $\Omega_1(e)$-positive Hermitian, and if we let 
	\begin{equation}\label{eq:heisenbergbracket}
	[(u,v),(u',v')]=(0,q_e(u,u'))\mathfa u,u'\in Z_{1/2}(e)\,,\,v,v'\in X_1(e)\ ,
	\end{equation}
	then $Z_{1/2}(e)\oplus X_1(e)$ is a Lie algebra isomorphic to $\ger h^e$ by the map $\phi^e$ from (\ref{eq:heisenbergiso}). 
	
	Since the subspaces $\ger g^e[\lambda]$, $\lambda=1,2$, are $\ger g^e[0]$-invariant, we obtain $\ger g^e[0]$-module structures on $Z_{1/2}(e)$ and $X_1(e)$ by transport of structure. Here, $\ger l^e=\ger g_0(e)\oplus\ger m^e$ centralises $\ger g^e[2]$, and $\ger g_1(e)$ acts on $X_1(e)$ via
	\[
		\delta.v=\delta(v)\ ,\ \xi_u^-.v=2(u\circ v)\mathfa u,v\in X_1(e)\,,\,\delta\in\ger k_1(e)\ .
	\]
	In particular, the action of $\ger g_1(e)$ on $\ger g^e[2]$ is equivalent to the action of $\ger{gl}(\Omega_1(e))$ on $X_1(e)$, and therefore faithful. 

 Futhermore, $\ger g^e[0]$ acts on $Z_{1/2}(e)$ via
	\begin{equation}\label{eq:g0g1action}
	\delta.v=\delta(v)\ ,\ \xi_u^-.v=-2\{uv^*e\}\mathfa\delta\in\ger k^e\,,\,u\in Z_0(e)\oplus X_1(e)\,,\,v\in Z_{1/2}(e)\ .
	\end{equation}
	In particular, $\ger z(\ger l^e\ltimes\ger h^e)=\ger z(\ger g_0(e)\ltimes\ger h^e)=\ger z(\ger h^e)=\ger g^e[2]$. 
\end{Prop}

\begin{proof}
The map $h_e$ is positive Hermitian by \cite[10.4]{loos_boundsymmdom}. Clearly, $[\ger h^e,\ger h^e]\subset\ger g^e[2]\subset\ger z(\ger h^e)$, and $\ger h^e$ is a generalised Heisenberg algebra. 
For $u,v\in Z_{1/2}(e)$, $[\eta_u^e,\eta_v^e]\in\ger g^e[2]$ and hence equals $\zeta_{-iw/2}^e$ \fs $w\in X_1(e)$. Since $\zeta_{-iw/2}^e(0)=\xi_{-iw/2}^-(0)=-\frac i2w$, 
\begin{align*}
-\tfrac i2w&=\Bracks0{\eta_u^e,\eta_v^e}(0)=[\xi_u^-,[\xi_e^-,\xi_v^-]](0)+[[\xi_e^-,\xi_u^-],\xi_v^-](0)\\
&=2\cdot\Parens0{\{ve^*u\}-\{ev^*u\}+\{eu^*v\}-\{ue^*v\}}=2\cdot\Parens0{\{vu^*e\}-\{uv^*e\}}\ .
\end{align*}
This proves that $Z_{1/2}(e)\oplus X_1(e)$ is a Lie algebra isomorphic to $\ger h^e$. 

For $x\in\Omega_1(e)^\circ$, let $b_x(u,v)=\Scp0{q_e(iu,v)}x$ \fa $u,v\in Z_{1/2}(e)$. Then $b_e$ is a symmetric bilinear form, positive definite since $\Omega_1(e)$ is regular and self-dual. Since $[(iu,0),(u,0)]=(0,q_e(iu,u))$ \fa $u\in Z_{1/2}(e)$, we find $\ger z(\ger h^e)=\ger g^e[2]=X_1(e)$.

Next, we consider the $\ger g^e[0]$-action on $X_1(e)$. If $\xi_u^-\in\ger p_0(e)\oplus\ger p_1(e)$ and $v\in X_1(e)$, then $[\xi_u^-,\zeta^e_{-iv/2}]=\zeta_{-iw/2}^e$ \fs $w\in X_1(e)$. We have 
\[
	-\tfrac i2w=\tfrac12[\xi_u^-,[\xi_e^-,\xi_{-iv/2}^-]]=\tfrac12\xi_{-i\trip0evu-i\trip0veu}(0)=-\tfrac i2\Parens0{\trip0evu+\trip0veu}\ .
\]

For $u\in Z_0(e)$, this is zero, and for $u\in X_1(e)$, it equals $-i(u\circ v)$ by (\ref{eq:tripvscircprod}). Similarly, for $\delta\in\ger k^e$, $[\delta,\zeta^e_{-iv/2}]=\zeta^e_{-iw/2}$ gives $-\tfrac i2w=[\delta,\xi_{-iv/2}^-](0)=-\tfrac i2\delta(v)$, so $w=\delta(v)$. For $\delta\in\ger k_0(e)$, this is zero, and so it is if $\delta\in\ger k^e$ is arbitrary and $v=e$. We have shown that $\ger g_0(e)$ centralises $\ger g^e[2]$, and that the $\ger g_1(e)$-action on $\ger g^e[2]$ is equivalent to the $\ger{gl}(\Omega_1(e))$-action on $X_1(e)$. In particular, $X_e^+$ is a cyclic vector of the $\ger g^e[0]$-module $\ger g^e[2]$. Since it is annihilated by $\ger m^e\subset\ger k^e$ and $\ger m^e$ is an ideal of $\ger g^e[0]$, we see that $\ger m^e$ centralises $\ger g^e[2]$. Evaluating $[\delta,\eta_v^e]$ and $[\xi_u^-,\eta_v^e]$ at zero \fa $\delta\in\ger k^e$, $u\in Z_0(e)\oplus X_1(e)$, and vectors $v\in Z_{1/2}(e)$ gives the remaining relations. 
\end{proof}

\begin{Lem}[jordanalg-auto-centre]
The centre of $\ger k_1(e)$ is trivial. In particular, if $\rk e\sge2$, then the derived algebra $\ger g_1(e)'=[\ger g_1(e),\ger g_1(e)]$ is a non-compact, non-Hermitian simple Lie algebra. If $\rk e\sle 1$, then $\ger g_1(e)=\reals\xi_e^-$ is Abelian.
\end{Lem}

\begin{proof}
By (\ref{eq:zerograding-decomp}), $\Ad(\gamma_e)(\ger g_1(e))=\gl0{\Omega_1(e)}$. The Lie algebra $\gl0{\Omega_1(e)}$ is reductive with centre $\reals M_e$. Because $X_1(e)$ is a simple Jordan algebra for $e\neq0$, $\gl0{\Omega_1(e)}'$ is a simple Lie algebra or zero. If $\rk e\sge2$, then there exists an idempotent $c\in X_1(e)$, $0<c<e$, and $M_c\subset\gl0{\Omega_1(e)}'$ generates an unbounded one-parameter group, so $\gl0{\Omega_1(e)}'$ is non-compact.

Finally, let $\delta\in\ger z(\ger k_1(e))$ and $u\in X_1(e)$. We have $0=[\delta,\jBox ue]=\jBox{(\delta u)}e=M_{\delta u}$, since $\delta e=0$, so $\delta u=0$. This shows that $\delta=0$.
\end{proof}

\paragraph{Principal faces}
Using the identification $\phi^e:Z_{1/2}(e)\oplus X_1(e)\to\ger h^e$ from \thmref{Prop}{heisenberg-alg}, we consider the cone $\Omega_1(e)\subset X_1(e)$ as a subset of $\ger g^e[2]=\ger z(\ger h^e)$. We point out that this notation is only meaningful if we keep the embedding $\phi^e$ attached to $e$ in mind. (Indeed, $\phi^e(-\Omega_1(e))$ and $\phi^{-e}(-\Omega_1(e))$ are distinct!) In what follows, the chosen embedding will always be clear from the context.

\begin{Prop}[nilpotent-cone]
We have $\Omega^\pm\cap\ger h^e=\Omega_1(e)$.
\end{Prop}

\begin{proof}
Let $\Omega$ be one of $\Omega^\pm\cap\ger h^e$. Then $\Omega$ is a closed pointed cone invariant under $N_G(\ger h^e)$, and in particular, under inner automorphisms of $\ger h^e$. Hence $\Omega\subset\ger z(\ger h^e)=\ger g^e[2]$ \cite[Lemma~I.13]{hno_conalheisenberg}, and by \thmref{Prop}{heisenberg-alg}, $\Omega$ is invariant under $\GL0{\Omega_1(e)}$. On the other hand, $X_e^+=\phi^e(e)\in\Omega$. Identifying $\Omega$ with its image in $X_1(e)$, this implies $\Omega_1(e)\subset\Omega$ and $\Omega^*\subset\Omega_1(e)^*=\Omega_1(e)$. Since $\Omega$ is pointed, the interior of $\Omega^*$ in $\ger g^e[2]$ is non-void. Hence, there is some $x\in\Omega^*\cap\Omega_1(e)^\circ$, and $\Omega_1(e)^\circ\subset\Omega^*$ since $\Omega_1(e)^\circ$ is homogeneous. It follows that $\Omega^*=\Omega_1(e)$, and by duality, $\Omega=\Omega_1(e)$. 
\end{proof}

\begin{Def}
Define $F_e^\pm=\Omega^\pm\cap(X_e^-)^\perp$. Since $-X_e^-\in\Omega^-\subset\Omega^+$, this is an exposed face of $\Omega^\pm$. We call $F_e^\pm$ a \emph{principal face}.
\end{Def}

\begin{Prop}[cone-parabolic]
We have $F_e^\pm=\Omega^\pm\cap\ger q^e$, and this is an exposed face of $\Omega^\pm$.
\end{Prop}

The \emph{proof} is preceded by two lemmata.

\begin{Lem}[nilpotent-face]
Let $e\sge c$ be non-zero tripotents, $n=\dim X_1(e)$, $k=\rk e$. Denote the canonical inner product of $X_1(e)$ by $\Scp0\cdot\cdot$. Then, \fa $u\in X_1(e)$, $v\in X_1(c)$, 
\[
	\Rscp0{\phi^e(u)}{\phi^c(v)}=\tfrac{4n}k\cdot\Scp0uv\nd\Rscp0{\phi^e(u)}{\phi^{-c}(v)}=0\ .
\]
\end{Lem}

\begin{proof}
Let $c\sle e$, $u\in X_1(e)$, $v\in X_1(c)$. Then $\vtheta(\zeta^e_{iu})=-\xi_{iu}^-+\tfrac12[\xi_e^-,\xi_{iu}^-]$, so
\begin{align*}
	\Rscp0{\phi^e(u)}{\phi^{\pm c}(v)}&=\Rscp0{\zeta^e_{-iu/2}}{\zeta^{\pm c}_{-iv/2}}\\
	&=-B(\xi_{iu/2}^--\tfrac12[\xi_e^-,\xi_{iu/2}^-],-\xi_{iv/2}^-\mp\tfrac12[\xi_c^-,\xi_{iv/2}^-])\\
	&=\tfrac14B(\xi_{iu}^-,\xi_{iv}^-)\mp\tfrac1{16}B(\xi_{iu}^-,[[\xi_c^-,\xi_{iv}^-],\xi_e^-])\ .\\
	\intertext{Since $\tfrac14[[\xi_c^-,\xi_{iv}^-],\xi_e^-]=-\xi_{iv}^-$ by (\ref{eq:fundid}), this is $0$ for $\phi^{-c}(v)$. For $\phi^c(v)$, by \thmref{Lem}{killing}, }
	&=\tfrac12B(\xi_{iu}^-,\xi_{iv}^-)=2\tr_Z(\jBox uv+\jBox vu)=\tfrac{4n}k\cdot\Scp0uv\ .
\end{align*}
\end{proof}

\begin{Lem}[gradingcone]
Let $\Omega\subset\ger g$ be a closed set invariant under $\reals_{\sge0}$ and $\Ad(\exp t\xi_e^-)$ \fa $t\in\reals$. If $\xi=\sum_{j=k}^\ell\xi_j\in\Omega$ where $\xi_j\in\ger g^e[j]$, then $\xi_k,\xi_\ell\in\Omega$. 
\end{Lem}

\begin{proof}
We have $\Ad(\exp t\xi_e^-)(\xi)=\sum_{j=k}^\ell e^{jt}\cdot\xi_j\in\Omega$ \fa $t\in\reals$. Hence, we have $\xi_k=\limk_{t\to\infty}e^{kt}\Ad(\exp -t\xi_e^-)(\xi)\in\Omega$ and $\xi_\ell=\limk_{t\to\infty}e^{-\ell t}\Ad(\exp t\xi_e^-)(\xi)\in\Omega$, proving the lemma.
\end{proof}

\begin{proof}[of \thmref{Prop}{cone-parabolic}]
If $e=0$, then $X_e^-=0$, $F_e^\pm=\Omega^\pm$, and $\ger q^e=\ger g$. W.l.o.g., we may assume $k=\rk e>0$. Since $\vtheta(\xi_e^-)=-\xi_e^-$, $\ad\xi_e^-$ is symmetric, and its eigenspaces are orthogonal. In particular, $\ger q^e\perp\ger g^e[-2]\ni X_e^-$, and $\Omega^\pm\cap\ger q^e\subset F^\pm_e$. 

For the converse, let $\xi\in F_e^\pm$, and write $\xi=\sum_{j=-2}^2\xi_j$ where $\xi_j\in\ger g^e[j]$. Since $X_e^-$ is an eigenvector of $\ad\xi_e^-$, $F_e^\pm$ is invariant under $\Ad(\exp t\xi_e^-)$ for all $t\in\reals$, so we can employ \thmref{Lem}{gradingcone}. In particular, $\xi_{-2}\in F_e^\pm\cap\ger h^{-e}$. By \thmref{Prop}{nilpotent-cone}, $\xi_{-2}=\phi^{-e}(u)$ for a unique $u\in-\Omega_1(e)$. By \thmref{Lem}{nilpotent-face},
\[
\Rscp0{\xi_{-2}}{X_e^-}=-\Rscp0{\phi^{-e}(u)}{\phi^{-e}(e)}=-\tfrac{2n}k\Scp0ue
\]
where $n=\dim X_1(e)$. This is positive if $u\neq0$, so $u=0$ and $\xi_{-2}=0$. Therefore, $\xi_{-1}\in\Omega^\pm\cap\ger h^{-e}$ by \thmref{Lem}{gradingcone}. Then \thmref{Prop}{nilpotent-cone} gives $\xi_{-1}=0$, and $\xi\in\ger q^e$. 
\end{proof}

\begin{Cor}[nilpotent-face]
	For any tripotent $e$, $\Omega_1(e)\subset\Omega^-$ is a face of $\Omega^+$ and $\Omega^-$. 
\end{Cor}

\begin{proof}
	It sufficient to show that $\Omega_1(e)$ is a face of $F_e^\pm$. Hence, let $\xi,\eta\in F_e^\pm$ \scth $\xi+\eta\in\Omega_1(e)$, and decompose $\xi=\sum_{j=0}^2\xi_j$, $\eta=\sum_{j=0}^2\eta_j$, according to the grading of $\ger q^e$. Then $\xi_0+\eta_0=0$ by assumption, and $\xi_0,\eta_0\in\Omega^\pm$ by \thmref{Lem}{gradingcone}. This implies $\xi_0=\eta_0=0$, and $\xi_1,\eta_1\in\Omega^\pm$ by the same lemma. But then \thmref{Prop}{nilpotent-cone} implies that $\xi_1=\eta_1=0$. Hence the claim. 
\end{proof}

\paragraph{Nilpotent faces and nilpotent orbits} 

We will now give a precise description of the conal nilpotent orbits. They are intimately related to the \emph{nilpotent faces} of $\Omega^\pm$. 

\begin{Def}
	Let $F\subset\Omega^\pm$ be a face. If $F^\circ$ contains a nilpotent (semi-simple) element of $\ger g$, we will call $F$ a \emph{nilpotent face} (\emph{semi-simple face}). 
	
	For any tripotent $e$, let $\mathcal O_e=\Ad(G)(X_e^+)$. Let $M_k$ be the set of rank $k$ tripotents.
\end{Def}

\begin{Th}[mincone-co-orbit]
Let $e$ be a tripotent of $\rk e=k$, and let $K^e=Z_K(e)$. Then 
\begin{equation}\label{eq:nilorbitface}
	\mathcal O_e=\textstyle\bigcup\nolimits_{c\in M_k}\Omega_1(c)^\circ=K\times_{K^e}\Omega_1(e)^\circ\ .
\end{equation}
In particular, $\mathcal O_e$ depends only on the rank of $e$; moreover,
\[
\rk e\neq\rk c\ \Rightarrow\ \mathcal O_e\cap\mathcal O_c=\vvoid\nd\rk e\sle\rk c\ \Rightarrow\ \mathcal O_e\subset\clos{\mathcal O_c}\ .
\]
Every nilpotent orbit in $\Omega^+$ is one of the $\mathcal O_c$; every nilpotent face is one of the $\Omega_1(c)$. 
\end{Th}

\begin{Rem}
The classification of conal nilpotent orbits is contained in \cite[Theorem~2]{vinberg}, \cite[Theorem~III.9]{hno_nilpotent}, \cite[Lemma~4]{paneitz_detcones}. The description in terms of tripotents and the connection to nilpotent faces is, however, new. While the parametrisation of the conal nilpotent orbits by tripotents might be deduced from \cite[Theorem~III.9]{hno_nilpotent} by applying \thmref{Prop}{discembed-tripotent}, our \emph{proof} of the more precise result is independent of existing results, and at the same time, shorter and more elementary.  
\end{Rem}

\begin{proof}[of \thmref{Th}{mincone-co-orbit}]
By \thmref{Prop}{sl2triple}, $\mathcal O_e\subset\Omega^-$. If $\rk e=\rk c$, then $\ell(e)=c$ \fs $\ell\in K$ \cite[Corollary~5.12]{loos_boundsymmdom}. Then $\Ad(\ell)(X_e^+)=X_c^+$ by (\ref{eq:kaction-nilpotent}), so $\mathcal O_e=\mathcal O_c$. 

Let $x\in\Omega^+$ be nilpotent, $x\neq0$. Then $x\in\Omega^-$ and there exists $g\in G$ such that $\Ad(g)(x)=x^+$ for some Cayley triple $(h,x^+,x^-)$ \cite[Theorems 9.2.1, 9.4.1]{collingwood-mcgovern}. By \thmref{Prop}{cayleyh1}, the triple is $(H_1)$, so $\Ad(g)(x)=X_e^+$ \fs tripotent $e$, by \thmref{Prop}{discembed-tripotent}. Let $F$ be the face of $\Omega=\Omega^\pm$ generated by $x$. Since $\Omega\cap\ger g^e[2]=\Omega_1(e)$ is a face of $\Omega$ by \thmref{Cor}{nilpotent-face}, $\Ad(g)(F)$ equals the face of $\Omega_1(e)$ generated by $x=X_e^+$. But this face is $\Omega_1(e)$ itself. 

By the Iwasawa decomposition, $G$ is generated by $K$ and the analytic subgroup $Q^e$ associated with $\ger q^e$. Since $Q^e$ normalises $\ger g^e[2]=X_1(e)$, we find that $\Ad(\ell)(F)=\Omega_1(e)$ \fs $\ell\in K$. From (\ref{eq:kaction-nilpotent}), $F=\Omega_1(c)$ \fs tripotent $c=\ell^{-1}(e)$ with $\rk c=k$. Let $G_1(c)$ be the analytic subgroup of $G$ associated with $\ger g_1(c)$. By \thmref{Prop}{heisenberg-alg}, the action of $G_1(c)$ on $\ger g^c[2]$ corresponds to the action of $\GL0{\Omega_1(c)}$ on $X_1(c)$, and is thus transitive on $F^\circ$. This proves the equation (\ref{eq:nilorbitface}), the exhaustion of nilpotent orbits in $\Omega^+$, and the exhaustion of nilpotent faces. Since $c$ is the only tripotent contained in $\Omega_1(c)^\circ$, $\Ad(G)(x)=\mathcal O_c$ does not contain any rank $k-1$ tripotent. Similarly, any tripotent $c'\sle c$ is contained in $\Omega_1(c)$, and therefore in $\overline{\mathcal O_c}$.  
\end{proof}

\begin{Cor}[minorbit-span]
	Let $e$ be primitive. Then $\Omega^-=\overline{\co0{\mathcal O_e}}=0\cup\co0{\mathcal O_e}$. 
\end{Cor}

\begin{proof}
	Let $C=\overline{\co0{\mathcal O_e}}\subset\Omega^-$. Then $C$ is a $G$-invariant closed convex cone. We have $\mathcal O_e=K\times_{K^e}(\reals_{>0}\cdot X_e^+)$, so $C=0\cup\co0{\mathcal O_e}=\reals_{\sge0}\cdot\co0{\Ad(K)(X_e^+)}$.  

To see that $C=\Omega^-$, it remains to be shown that $iz\,\frac\partial{\partial z}\in C$. We have $\pm X^\pm_e\in\mathcal O_e$, so $i\jBox ee=\tfrac12(X_e^+-X_e^-)\in C$. By \thmref{Lem}{extremalgens} and the $K$-invariance of $C$, $\omega^-\subset C$. But $iz\,\frac\partial{\partial z}\in\omega^-$, and therefore, $C=\Omega^-$.
\end{proof}

\begin{Cor}
	Every conal nilpotent orbit $\mathcal O_e$ is a $K$-equivariant fibre bundle over $K/K^e=M_k$ (for $k=\rk e$), with contractible fibres. In particular, $\mathcal O_e$ is $K$-equivariantly homotopy equivalent to $M_k$. 
\end{Cor}

\begin{Rem}
	The projection of the fibre bundle $\Omega_1(e)^\circ\to\mathcal O_e\to M_k$ ($k=\rk e$) associates to a nilpotent $x$ the \emph{unique} $y$ which generates the same face of $\Omega^-$ as $x$ and is the nilpositive element of a Cayley triple. In particular, with any nilpotent element of $\Omega^-$, we may associate a \emph{canonical} Cayley triple.  
\end{Rem}

\section{Classification of the faces of the minimal invariant cone}

In this section, we classify all faces of $\Omega^-$. First, we study $F_e^\pm$ in detail. 

\subsection{Fine structure of the principal faces} 

We have seen that the exposed face $F_e^\pm$ is contained in the maximal parabolic $\ger q^e$, and in particular, invariant under inner automorphisms of $\ger q^e$. However, this is not the definitive statement on $F_e^\pm$: the linear span of $F_e^\pm$ is a proper ideal of $\ger q^e$. 

\begin{Prop}[cone-gradingzero]
We have $F_e^\pm=\Omega^\pm\cap(\ger g_0(e)\ltimes\ger h^e)$. If $\rk e<r$, then both of the faces $F_e^\pm$ span $\ger g_0(e)\ltimes\ger h^e$. If $\rk e=r$, then $F_e^\pm=\Omega^\pm\cap\ger g^e[2]=\Omega_1(e)\subset\Omega^-$. 
\end{Prop}

The \emph{proof} requires a preparatory lemma. Fix a frame $e_1,\dotsc,e_r$, and recall the compact Cartan subalgebra $\ger t=\ger t^+\oplus\ger t^-$ from (\ref{eq:cptcsa}). Let $\ger a=\Span0{\xi_{e_1}^-,\dotsc,\xi_{e_r}^-}_\reals$, $\ger m=\ger z_{\ger k}(\ger a)$. 

\begin{Lem}[gradingzero-cartan]
Let $e=\mathbf e_k=e_1+\dotsm+e_k$. We have $\ger t\cap\ger m=\ger t^+$, and
\[
\ger t^e[0]=\ger t\cap\ger g^e[0]=\langle i\cdot e_j^{\phantom*}\Box\,e_j^*\mid j=k+1,\dotsc,r\rangle\oplus\ger t^+\subset\ger l^e\ .
\]
The subalgebras $\ger g^e[0]$, $\ger g_0(e)$, and $\ger m^e$ of $\ger g$ are $\ger t$-invariant. Moreover, $\ger t_0(e)=\ger t\cap\ger g_0(e)$ and $\ger t^+\cap\ger m^e$ are Cartan subalgebras of $\ger g_0(e)$ and $\ger m^e$, respectively. 
\end{Lem}

\begin{proof}
Since $[\delta,\xi_{e_j}^-]=\xi_{\delta e_j}^-$ \fa $\delta\in\ger k$, $\ger t^+\subset\ger m=\ger z_{\ger k}(\ger a)$. For the converse, we have $\trip0ccc=c\neq0$ for $c=e_j$, so $\ger m\cap\ger t^-=0$ and $\ger t\cap\ger m=\ger t^+$. Since $\ger m^e\subset\ger m$, $\ger t\cap\ger m^e\subset\ger t^+$. 

Moreover, $ie_j\,\Box\,e_j^*\in\ger k_0(e)$ if $j>k$, and if $j\sle k$, then 
\[
[\delta,ie_j\,\Box\,e_j^*]=i\cdot(\delta e_j)\,\Box\,e_j^*+i\cdot e_j\,\Box\,(\delta e_j)^*=0\mathfa\delta\in\ger k^e\ ,
\]
and $[\xi_u^-,i\cdot e_j\,\Box\,e_j^*]=-\xi_{i\{e_j^{\phantom*}e_j^*u\}}^-=0$ \fa $u\in Z_0(e)$.

We conclude that $\ger l^e$ is $\ger t$-invariant, and $\ger t^e[0]=\langle ie_j^{\phantom*}\Box\,e_j^*\mid j=k+1,\dotsc,r\rangle\oplus\ger t^+\subset\ger l^e$. In addition, $\ger t_0(e)$ and $\ger t^+\cap\ger m^e$ are Cartan subalgebras of $\ger g_0(e)$ and $\ger m^e\;$, respectively \cite[Chapter~VIII, \S~3.1, Proposition~3]{bourbaki}.
\end{proof}

Let $\Omega_0^\pm(e)$ denote the minimal and maximal cones of the Lie algebra $\ger g_0(e)$, cf.~\thmref{Def}{facial}. Likewise, set $\omega_0^\pm(e)=\Omega_0^\pm(e)\cap\ger t_0(e)$. Then 
\[
\omega_0^+(e)=\omega_0^-(e)^*\nd\omega_0^-(e)=\langle iH_\alpha\mid\alpha\in\Delta_n^{++}\,,\,\ger g^\alpha_\cplxs\subset\ger g_{0,\cplxs}(e)\rangle\ .
\]
Here, $\ger g_{0,\cplxs}(e)=\ger g_0(e)\otimes\cplxs$. The set $\{\alpha\in\Delta_n^{++}\mid\ger g^\alpha_\cplxs\subset\ger g_0(e)\}$ coincides with the set of positive non-compact roots for $\ger g_0(e)$, since this algebra is $\ger t$- and $\vtheta$-invariant \cite[Chapter~VIII, \S~3.1, Proposition~3]{bourbaki}. 

\begin{proof}[of \thmref{Prop}{cone-gradingzero}]
We have $\rk e<r$ if and only if $\ger g_0(e)\neq0$. In this case, $\ger h=\ger t_0(e)\oplus X_1(e)$ is a compact Cartan subalgebra of $\ger g_0(e)\ltimes\ger h^e$. The intersection of a generating cone with such a Cartan subalgebra completely determines the cone \cite[Proposition~III.5.14 (ii)]{hilgert-hofmann-lawson}. Thus, we claim that $F_e^\pm=\Omega^\pm\cap(\ger g_0(e)\ltimes\ger h^e)$, independent of the rank of $e$. This will imply the assertion for $\rk e<r$; for $\rk e=r$, it follows from \thmref{Prop}{nilpotent-cone}. 

Assume that we have shown $F_e^\pm\subset\ger l=\ger g_0(e)\ltimes\ger h^e$ and that $F_e^\pm\cap\ger h$ is solid in $\ger h$. Since $F_e^\pm=\Omega^\pm\cap\ger q^e$, $F_e^\pm$ is invariant under inner automorphisms of $\ger q^e$, and in particular, of $\ger l$. It follows that $F_e^\pm$ is the unique pointed invariant cone in $\ger l$ whose intersection with $\ger h$ is $F_e^\pm\cap\ger h$, and this intersection is regular in $\ger l$ \cite[Theorem~III.5.15, Proposition III.5.14 (iii)]{hilgert-hofmann-lawson}. Thus, once we have shown our assumption, it is clear that $\ger l$ is spanned by $F_e^\pm$. 

In view of \thmref{Lem}{gradingcone}, it is sufficient to prove that $\Omega^\pm\cap\ger g^e[0]=\Omega_0^\pm(e)$, and that $\omega^\pm\cap\ger g^e[0]=\omega_0^\pm(e)$. Moreover, we may assume $e=\mathbf e_k=e_1+\dotsm+e_k$, and $k\sge1$. From (\ref{eq:zerograding-decomp}), we have $\ger g^e[0]=\ger g_0(e)\oplus\ger m^e\oplus\ger g_1(e)$ \fs compact reductive ideal $\ger m^e\subset\ger m=\ger z_{\ger k}(\ger a)$ of $\ger g_0(e)$. Moreover, $\ger g_0(e)\oplus\ger m^e$ is invariant under $\ger t$ by \thmref{Lem}{gradingzero-cartan}. Let $p_{\ger t}$ be the orthogonal projection onto $\ger t$. Since $\ger k\perp\ger p$ and 
\[
\Rscp0{i\cdot e_j^{\phantom*}\,\Box\,e_j^*}\delta=-2i\tr_Z((\delta e_j^{\phantom*})\,\Box\,e_j^*)=0\mathfa\delta\in\ger k\cap\ger g^e[0]\,,\,j\sle k\ ,
\]
by \thmref{Lem}{killing}, $p_{\ger t}$ leaves $\ger g^e[0]$ invariant. Thus, $p_{\ger t}(\Omega^-\cap\ger g^e[0])=\omega^-\cap\ger g^e[0]\subset\ger t^e[0]$. 

\thmsref{Lem}{extremalgens}{gradingzero-cartan} give $\omega^-\cap\ger g_0(e)=\omega_0^-(e)$. Hence, $\Omega^-\cap\ger g_0(e)=\Omega_0^-(e)$ \cite[Theorem~2]{paneitz_detcones}. Let $\widetilde\Omega^\pm=\Omega^\pm\cap\ger g^e[0]$. Then $\widetilde\Omega^\pm$ is closed, pointed, and invariant under inner automorphisms. Hence, $\ger a^\pm=\widetilde\Omega^\pm-\widetilde\Omega^\pm$ is an ideal of $\ger g^e[0]$. Since $\ger g^e[0]$ is reductive, so is $\ger a^\pm$, and moreover, quasihermitian \cite[Proposition~II.2 and Lemma~II.4]{neeb_invconv}. \thmref{Lem}{jordanalg-auto-centre} implies $\ger a^\pm\cap\ger g_1(e)=0$, since $\ger a^\pm$ has neither proper non-compact Abelian nor non-Hermitian simple ideals. We conclude $\widetilde\Omega^\pm\subset\ger a^\pm\subset\ger g_0(e)\oplus\ger m^e$.

Let $\xi\in\omega^+\cap\ger m^e$. Seeking a contradiction, assume $\xi\neq0$. Then there is $\alpha\in\Delta_n^{++}$ such that $\alpha(\xi)>0$. Since $[\xi_e^-,\xi]=0$, $\ger g^\alpha_\cplxs\subset[\xi,\ger p^+]\subset\ger p^+$ is $\ad\xi_e^-$-invariant, and hence contained $\ger g^e[\ell]_\cplxs$ \fs $\ell$. But $[\xi_e^-,\ger g_\cplxs^\alpha]\subset[\ger p_\cplxs,\ger p_\cplxs]\subset\ger k_\cplxs$, and $\ger k_\cplxs\cap\ger g_\cplxs^\alpha=0$, so necessarily $\ell=0$. Since $\ger m^e$ is an ideal of $\ger g^e[0]$, so $\ger g^\alpha_\cplxs=[\xi,\ger g_\cplxs^\alpha]\subset\ger m^e\cap\ger p^+=0$, a contradiction. Therefore, $\omega^+\cap\ger m^e=0$. 

Since $\ger t\cap(\ger g_0(e)\oplus\ger m^e)=\ger t_0(e)\oplus\ger m^e\cap\ger t^+$, the projections $p_{\ger t}$ and $p_{\ger m^e}$ commute, and $p_{\ger m^e}(\omega^\pm\cap\ger g^e[0])=\omega^\pm\cap\ger m^e=0$. Consequently, $\omega^-\cap\ger g^e[0]=\omega_0^-(e)$, and this entails $\widetilde\Omega^-=\Omega_0^-(e)$. As for the dual cone, clearly $\Omega^+\cap\ger g^e[0]\subset\Omega_0^-(e)^*=\Omega_0^+(e)$. In particular, we have the inclusion $\omega^+\cap\ger g^e[0]\subset\omega_0^+(e)$. 

Conversely, for $\alpha\in\Delta_n^{++}$, non-vanishing on $\omega_0^+(e)$, we have $\ger g^\alpha_\cplxs\subset\ger p^+\cap\ger g^e[0]_\cplxs$. If $\ger g^\alpha_\cplxs\not\subset\ger g_0(e)_\cplxs$, then, since $\ger g^\alpha_\cplxs\cap\ger m^e=0$, $\ger g^\alpha_\cplxs\subset\ger g_1(e)_\cplxs$. Because $\alpha(\ger t_0(e))\neq0$, we find that 
$\ger g^\alpha_\cplxs\subset[\ger t_0(e),\ger g_1(e)_\cplxs]\subset[\ger g_0(e)_\cplxs,\ger g_1(e)_\cplxs]=0$, which is a contradiction. Therefore, $\ger g^\alpha_\cplxs\subset\ger g_0(e)_\cplxs$. This means that $\alpha$ is a root for $\ger g_0(e):\ger t_0(e)$, and hence $-i\alpha\sge0$ on $\omega_0^+(e)$ by definition. We have established that $-i\alpha\sge0$ on $\omega_0^+(e)$, \fa $\alpha\in\Delta_n^{++}$. Hence, we have that $\omega_0^+(e)\subset\omega^+\cap\ger g^e[0]$, and equality follows. In particular, we have $\Omega^+\cap\ger g^e[0]=\Omega_0^+(e)$ \cite[Theorem~2]{paneitz_detcones}.
\end{proof}

\begin{Rem}
We use techniques due to Neeb \cite[Proposition~VIII.3.30]{neeb_holconv}.
\end{Rem}

\subsection{Semi-simple and general faces}

We now construct the semi-simple faces and use general results on Lie algebras with invariant cones to determine the structure of arbitrary faces of $\Omega^\pm$. In particular, all these faces span subalgebras of $\ger g$ whose Levi complements are given by the $\ger g_0(e)$. 

\paragraph{Construction of semi-simple faces}

\begin{Prop}
	We have $\Omega_0^\pm(e)=F_e^\pm\cap\Omega_1(e)^\perp=\Omega^\pm\cap(X_e^-)^\perp\cap(X_e^+)^\perp$, and this set is an exposed semi-simple face of $\Omega^\pm$. 
\end{Prop}

\begin{proof}
	We have $\Omega_1(e)\subset\Omega^-\subset\Omega^+$, so that $\Omega^\pm\cap\Omega_1(e)^\perp$ is an exposed face of $\Omega^\pm$. We have $\Omega^\pm\cap\Omega_1(e)^\perp=\Omega^\pm\cap(X_e^+)^\perp$, because $X_e^+=\phi^e(e)\in\Omega_1(e)^\circ$. 
	
	As the intersection of exposed faces, $F=F_e^\pm\cap\Omega_1(e)^\perp$ is exposed. Since $\Omega_1(e)$ spans $\ger g^e[2]$, \thmref{Lem}{gradingcone} and \thmsref{Prop}{nilpotent-cone}{cone-gradingzero} show that $F=\Omega_0^\pm(e)$. 
\end{proof}

\begin{Cor}
	The nilpotent faces of $\Omega^\pm$ are exposed. 
\end{Cor}

\begin{proof}
	Note $\Omega_1(e)=F_e^\pm\cap\Omega_0^\mp(e)^\perp$, and exposed faces form a complete lattice.
\end{proof}

We will show that the $\Omega_0^\pm(e)$ exhaust the set of semi-simple faces. In view of the following lemma, it will suffice to show that they exhaust them up to conjugacy.

\begin{Lem}[g0conj]
	Let $\ger h$ be a subalgebra of $\ger g$ conjugate to $\ger g_0(e)$ \fs tripotent $e$. Then there exists a tripotent $c$, $\rk c=\rk e$, \scth $\ger h=\ger g_0(c)$. 
\end{Lem}

\begin{proof}
	We may assume that $\ger h\neq0$. Recall that $\ger g$ is the set of all complete holomorphic vector fields on the bounded symmetric domain $D\subset Z$. The group $G$ acts on the set of faces of $D$, and each of the faces is of the form $F=e+D_0(e)$ where $e$ is the unique tripotent contained in $F$. The normaliser of the face $F$ is the parabolic $\ger q^e$, and the latter is invariant under $\ad\xi_e^-$. 
	
	The unique $\ad\xi_e^-$-invariant complement of the nilradical of $\ger q^e$ is $\ger g^e[0]$, and $\ger g_0(e)$ is the unique Hermitian simple ideal therein by (\ref{eq:zerograding-decomp}) and \thmref{Lem}{jordanalg-auto-centre}. By assumption, $\ger h$ is the Hermitian simple ideal in the canonical complement of the nilradical of the normaliser of a face of $D$.\end{proof}
	
\paragraph{The structure of general faces}

\begin{Lem}[face-subalg]
	Let $H$ be a Lie group, and $\Omega\subset\ger h$ a closed convex $\Ad H$-invariant cone. Any face $F$ of $\Omega$ spans a subalgebra of $\ger h$. In fact, if $\xi\in F^\circ$ and $\eta\in\ger n_{\ger h}(\reals\xi)$, then $\ad\eta$ leaves $\Span0F_\reals$ invariant. 
\end{Lem}

\begin{proof}
	Let $F\subset\Omega$ be a face, and $\xi\in F^\circ$. Let $\eta\in\ger n_{\ger h}(\reals\xi)$. Then \fa $t$, $\Ad(\exp t\eta)$ normalises $\reals\xi$. Furthermore, $G=\Ad(t\exp\xi)(F)$ is a face of $\Omega$, since $\Ad(t\exp\eta)$ is a linear automorphism of $\ger h$ leaving $\Omega$ invariant. Moreover, $\Ad(t\exp\eta)$ is an open map, so $\xi\in G^\circ$. Hence, $G=F$ \cite[Theorem~13.1]{rockafellar}, and differentiating with respect to $t$, we obtain $[\eta,F]\subset \clos{\reals_{>0}\cdot F-F}=\Span0F_\reals$. In particular, we may choose $\eta=\xi$. Since $F^\circ$ is dense in $F$, the claim follows. 
\end{proof}

\begin{Def}
	Let $F$ be a face of $\Omega^\pm$. We let $\ger g_F$ be the subalgebra spanned by $F$ and call this the \emph{face algebra}. Furthermore, let $\ger r_F$ be the radical of $\ger g_F$, $\ger n_F$ the nilradical, $\ger z_F=\ger z(\ger g_F)$ the centre, and let $G_F$ be the analytic subgroup of $G$ associated with $\ger g_F$.  
\end{Def}

\begin{Prop}[face-structure]
	Let $\Omega=\Omega^\pm$ and $F$ be a face of $\Omega$. There exists a compact Cartan subalgebra $\ger t_F\subset\ger g_F$, and a unique $\ger t_F$-invariant Levi complement $\ger s_F$. Then $\ger s_F$ is quasihermitian semi-simple, and $\ger g_F=\ger s_F\ltimes\ger n_F$ where $[\ger n_F,\ger n_F]\subset\ger z_F$. Moreover, $\ger n_F=[\ger t_F\cap\ger s_F,\ger n_F]\oplus\ger z_F$ as vector spaces. If $\ger z_F=0$, then $\ger g_F=\ger s_F$; if $\ger g_F$ is solvable, then $\ger g_F=\ger z_F$ is Abelian.
\end{Prop}

If the \emph{proof}, we will need the following definition.

\begin{Def}
	Let $\ger a$ be a real Lie algebra with compactly embedded Cartan subalgebra $\ger b$. Let ${}^*$ be the complex conjugation of $\ger a_\cplxs$ with respect to $i\ger a$. A root $\alpha$ of $\ger a:\ger b$ is called \emph{compact} if $\alpha([x,x^*])>0$ \fs $x\in\ger a_\cplxs^\alpha$, and \emph{non-compact} otherwise. Moreover, $\ger a$ said to have \emph{cone potential} if $[x,x^*]\neq0$ for each non-zero non-compact root vector $x$ \cite[Definition~VII.2.22]{neeb_holconv}.
\end{Def}

\begin{proof}[of \thmref{Prop}{face-structure}]
The face $F$ is an $\Ad G_F$-invariant closed regular convex cone in $\ger g_F$. It follows that $\ger g_F$ is quasihermitian with a compactly embedded Cartan subalgebra $\ger t_F$, and a maximal compactly embedded subalgebra $\ger k_F$.

Let $\ger r_F$ be the radical of $\ger g_F$. There exists a unique $\ger t_F$-invariant Levi complement $\ger s_F$ of $\ger r_F$ \cite[Propositions VII.1.9, VII.2.5]{neeb_holconv}, and it is also $\ger k_F$-invariant. Furthermore, we have $\ger k_F=\ger r_F\cap\ger k_F\oplus\ger s_F\cap\ger k_F$, and if $\ger l_F=\ger s_F\oplus\ger t_F\cap\ger r_F\cap\ger z(\ger g_F)^\perp$, then $\ger l_F$ is a reductive subalgebra which is complementary to $\ger n_F$ in $\ger g_F$ (\emph{loc. cit.}). This subalgebra is quasihermitian \cite[Lemma~VIII.3.5, Theorem VIII.3.6]{neeb_holconv}, and hence, the sum of a compact and of Hermitian simple ideals. 

Since $F$ is an invariant regular cone in $\ger g_F$, this Lie algebra has cone potential \cite[Theorem~III.6.18]{hilgert-hofmann-lawson}. Then $[\ger n_F,\ger n_F]\subset\ger z_F=\ger z(\ger g_F)$ \cite[Theorem~III.6.23]{hilgert-hofmann-lawson}. The subset $\ger l_F\cap\ger t_F\oplus\ger z_F$ is a Cartan subalgebra of $\ger g_F$ \cite[Theorem~VII.2.26]{neeb_holconv}, and since it is contained in $\ger t_F$, we find $\ger t_F=\ger l_F\cap\ger t_F\oplus\ger z_F$. 

Let $\Delta_F=\Delta(\ger g_{F,\cplxs}:\ger t_{F,\cplxs})$, and denote by $\Delta_{F,n}$ and $\Delta_{F,c}$ the subsets of non-compact and compact roots, respectively. There exists a unique adapted positive system $\Delta_F^{++}$ \scth $\omega_F^-\subset F\cap\ger t_F\subset\omega_F^+$ where $\omega_F^-$ is the cone spanned by $i[x^*,x]$ for $x\in\ger g_{F,\cplxs}^\alpha$, $\alpha\in\Delta_{F,n}^{++}$, and $\omega_F^+$ is the set of all $H\in\ger t$ \scth $-i\alpha(H)\sge0$, \fa $\alpha\in\Delta_{F,n}^{++}$ \cite[Theorem~VII.3.8]{neeb_holconv}. Let $p_{\ger t,F}:\ger g_F\to\ger t_F$ be the projection along $[\ger t_F,\ger g_F]$. Then we have $F=\Set1{\xi\in\ger g_F}{p_{\ger t,F}(\Ad(G_F)(\xi))\subset\ger t_F\cap F}$ (\emph{loc.~cit.}). 

Thus, $\Omega_F^-=\Set1{\xi\in\ger g_F}{p_{\ger t,F}(\Ad(G_F)(\xi))\subset\omega_F^-}$ is contained in $F$, and hence, pointed. It follows that $\Omega_F^-$ is the closed convex hull of $\Ad(G_F)(\omega_F^-)$ \cite[Corollary~VIII.3.31]{neeb_holconv}, so it is solid in $\ger g_F$ \cite[Proposition~III.5.14]{hilgert-hofmann-lawson}. Consequently, $\ger l_F$ has no non-zero compact ideal \cite[Proposition~III.3.30]{neeb_holconv}, and we have $\ger l_F=\ger s_F$. Hence, $\ger g_F=\ger s_F\oplus[\ger t_F\cap\ger s_F,\ger n_F]\oplus\ger z_F$ as vector spaces, because $\ger t_F=\ger s_F\cap\ger t_F\oplus\ger z_F$. Since any Abelian ideal of $\ger g_F$ is central \cite[Proposition~VII.3.23]{neeb_holconv}, we deduce that $\ger g_F=\ger s_F$ if $\ger z_F=0$. By the same token, $\ger g_F$ is Abelian if it is solvable. 
\end{proof}

Next, we determine the structure of the Levi complement $\ger s_F$. 

\begin{Prop}[face-semisimplepart]
	Let $\Omega=\Omega^\pm$ and $F\subset\Omega$ be a face. Let $\ger t_F$ be a compact Cartan subalgebra of $\ger g_F$ and $\ger s_F$ denote the $\ger t_F$-invariant Levi complement of $\ger g_F$. There exists a tripotent $e$ \scth $\ger s_F=\ger g_0(e)$. 
\end{Prop}

The \emph{proof} requires a little spadework. We begin with three lemmata which reduce the question to the study of the extremal rays of the cone $\ger s_F\cap\omega^-$. 

\begin{Lem}[sl2triple-face]
	Let $\ger s$ be an ideal of $\ger s_F$, and $i\cdot\jBox ee\in\ger s$. Then $\ger s^e\subset\ger s$ (where we recall $\ger s^e=\Span0{\xi_e^-,X_e^\pm}_\reals$ from \thmref{Prop}{sl2triple}).  
\end{Lem}

\begin{proof}
	Since $i\cdot\jBox ee=\tfrac12(X_e^+-X_e^-)$ and $\pm X_e^\pm\in\Omega^-$, we have $i\cdot\jBox ee\in\Omega^-\cap\ger s\subset F$. Because $F$ is a face, it follows that $\pm X_e^\pm\in F$. Now, $i\cdot\jBox ee$ remains unchanged if we replace $e$ by $te$ where $t\bar t=1$. \thmref{Th}{mincone-co-orbit} shows that the minimal nilpotent orbit of $\Omega^-\cap\ger s^e$ is the union of the rays spanned by the $X_{te}^-$, $t\bar t=1$. By \thmref{Cor}{minorbit-span}, the minimal cone $\Omega^-\cap\ger s^e$ is generated by this orbit, and hence $\ger s^e\subset\ger g_F$. Now, $\ger s^e$ it is not completely contained in $\ger n_F$ and is simple, so it is contained in $\ger s_F$. Since it intersects $\ger s$ non-trivially, we conclude $\ger s^e\subset\ger s_F$. 
\end{proof}

\begin{Lem}[smallertriple-face]
	Let $\ger s$ be an ideal of $\ger s_F$ and $i\cdot\jBox ee\in\ger s$. Then $\ger s^c\subset\ger s$ \fa $0<c\sle e$.
\end{Lem}

\begin{proof}
	From the previous lemma, we have $\pm X_e^\pm\in F\cap\ger s_F$. Then $F$ contains the faces $\Omega_1(\pm e)$ of $\Omega$ generated by these vectors. In particular, $\pm X_c^\pm\in F$ \fa $0<c\sle e$, and $i\cdot\jBox cc=\tfrac12(X_c^+-X_c^-)\in F$. The simple algebra $\ger s^c$ cannot be completely contained in $\ger n_F$. Arguing as in the previous lemma, we find $\ger s^c\subset\ger s_F$. 
\end{proof}

\begin{Lem}[face-simple]
	Assume that the span of those $i\cdot\jBox ee$ which belong to $\ger s_F$ contains a Cartan subalgebra of $\ger s_F$. Then the algebra $\ger s_F$ is simple.
\end{Lem}

\begin{proof}
	Assume that $\ger s^F$ splits as the direct sum of ideals $\ger s_1\oplus\ger s_2$. By assumption, there exist orthogonal tripotents $e_j$ \scth $i\cdot\jBox{e_j}{e_j}\in\ger s_j$. But then $e=e_1+e_2$ satisfies $\jBox ee=\jBox{e_1}{e_1}+\jBox{e_2}{e_2}$. By \thmref{Lem}{sl2triple-face}, $\ger s^e\subset\ger s_F$. Since $\ger s^e$ is simple, it must be contained in one of the ideals, $\ger s^e\subset\ger s_1$ (say). But then $i\cdot\jBox{e_2}{e_2}\in\ger s_1$, by \thmref{Lem}{smallertriple-face}, a contradiction!
\end{proof}

\begin{proof}[of \thmref{Prop}{face-semisimplepart}]
The semi-simple subalgebra $\ger s=\ger s_F$ is reductive in $\ger g$ and, possibly replacing $F$ by a $G$-conjugate, we may assume that it is $\vtheta$-invariant \cite[Lemma 1.1.5.5]{warner_semsimlg}. Then we have $\ger k_F\cap\ger s\subset\ger k$ \cite[Proposition~VII.2.5]{neeb_holconv}. Replacing $\ger t_F$ by a conjugate under inner automorphisms of $\ger k_F\cap\ger s$ (which are elements of $K$), we may assume $\ger t\cap\ger s\subset\ger t_F$. Then $\ger t_F\cap\ger s$ is contained in a Cartan subalgebra of $\ger g$ contained in $\ger k$. Replacing $F$ by a $K$-conjugate, we may assume that $\ger t_F\cap\ger s\subset\ger t$, so $\ger t\cap\ger s=\ger t_F\cap\ger s$ is a Cartan subalgebra of $\ger s$ contained in $\ger k_F\cap\ger s\subset\ger k$. 

It follows that $\Delta_{\ger s}=\Delta(\ger s_\cplxs:\ger s_\cplxs\cap\ger t_\cplxs)\subset\Delta$, and that the subsets $\Delta_{\ger s,c}$ and $\Delta_{\ger s,n}$ of compact and non-compact roots are, respectively, contained in $\Delta_c$ and $\Delta_n$. We may choose an adapted positive system $\Delta_{\ger s}^{++}$ contained in $\Delta^{++}$. Let $\omega_{\ger s}^-$ be the cone spanned by $iH_\alpha$ where $\alpha\in\Delta_{\ger s,n}^{++}$, and let $\omega_{\ger s}^+$ be the set of all $H\in\ger t\cap\ger s$ \scth $-i\alpha(H)\sge0$ \fa $\alpha\in\Delta_{\ger s,n}^{++}$. It is immediate that $\Delta_{\ger s,n}^{++}$ is the set of all $\alpha\in\Delta_n^{++}$ for which $\alpha(\ger t\cap\ger s)\neq0$, and hence $\ger s\cap\omega^\pm=\omega_{\ger s}^\pm$. Since $\omega^-$ is pointed, $\omega_{\ger s}^-$ is pointed, and its dual cone $\omega_{\ger s}^+$ is solid in $\ger t\cap\ger s$. Since $\omega^+$ is pointed, $\omega_{\ger s}^+$ is pointed, and its dual cone $\omega_{\ger s}^-$ is solid in $\ger t\cap\ger s$. Therefore, both of $\omega_{\ger s}^\pm$ are regular cones in $\ger t\cap\ger s$. 

Since $\Omega$ is an invariant regular cone, $\ger g$ has cone potential \cite[Theorem~III.6.18]{hilgert-hofmann-lawson}. We have $\Delta_{\ger s,n}\subset\Delta_n$, so $\ger s$ has cone potential, too. Since $\ger s$ is semi-simple, there exist unique invariant convex cones $\Omega^\pm_{\ger s}\subset\ger s$ \scth $\Omega^\pm_{\ger s}\cap\ger t=\omega_{\ger s}^\pm$ \cite[Theorem~III.5.15, Proposition III.5.14]{hilgert-hofmann-lawson}, and they are regular. Since $\Omega\cap\ger t\cap\ger s=\omega_{\ger s}^\pm$, it follows that $\Omega\cap\ger s=\Omega_{\ger s}^\pm$ is an invariant regular convex cone in $\ger s$. Because $\Omega_{\ger s}^-$ is pointed, $\ger s$ has no compact ideals, and is therefore a Hermitan non-compact Lie algebra \cite[Proposition~VIII.3.30]{neeb_holconv}. 

Observe now that $\ger s$ is $\vtheta$-invariant, and that $\ger s_\cplxs\cap\ger p_\cplxs=\ger s_\cplxs\cap\ger p^+\oplus\ger s_\cplxs\cap\ger p^-$ because $\Delta_{\ger s,n}\subset\Delta_n$. This decomposition allows the reconstruction of the triple product, and it follows that $Z_F=\ger s_F\cap\ger p$, which is a positive Hermitian Jordan triple in its own right, $\ger s_F$ being Hermitian non-compact, is a subtriple of $\ger p=Z$. Because $\omega_{\ger s}^+=\ger s\cap\omega^+$, it follows from \thmref{Cor}{extremalgens-tripotents} that $\ger t\cap s$ is spanned by those $i\cdot\jBox ee$ which lie in $\ger t\cap\ger s$. By \thmref{Lem}{face-simple}, $\ger s$ is simple, so $Z_F$ is simple.  

After renumbering, $ie_j^{\phantom*}\Box\,e_j^*\in\ger t\cap\ger s$ for $j=1,\dotsc,r_F$ where $e_1,\dotsc,e_{r_F}$ forms a frame of $Z_F$. If $e=e_{j+1}+\dotsm+e_r$, then the simplicity of $Z_F$ implies $Z_F=Z_0(e)$. Thus, $\ger s=\ger g_0(e)$. Finally, invoking \thmref{Lem}{g0conj}, this conclusion also holds without replacing $F$ by a $G$-conjugate.
\end{proof}

\subsection{Determination of the faces with non-reductive face algebra}

In order to determine all faces with non-reductive face algebra, the main step is to understand their centres. This is the content of the following proposition, which also will help us determine the faces with reductive face algebra. 

\begin{Prop}[facestructure2]
	Let $F\subset\Omega=\Omega^\pm$ be a face and $\ger g_F=\ger g_0(e)\ltimes\ger n_F$ where $\rk e<r$. Assume that $\ger z_F=\ger z_{\ger g_F}(\ger g_0(e))$. Possibly replacing $e$ by $-e$, we have $\ger g_F\subset\ger q^e$, $\ger n_F\subset\ger h^e$, and there exists a unique $c\sle e$ \scth $\ger z_F=\ger g^c[2]$ and $F\cap\ger z_F=\Omega_1(c)$. 
\end{Prop}

The \emph{proof}Êrequires some preparatory lemmata. 

\begin{Lem}[commjordanprod]
	Let $u,v\in X_1(e)$. Then $[\phi^e(u),\phi^{-e}(v)](0)=u\circ v$ where $\circ$ is the Jordan algebra product of $X_1(e)$. 
\end{Lem}

\begin{proof}
	Recall from (\ref{eq:heisenbergiso}) that $\phi^{\pm e}(u)=\zeta_{-iu/2}^{\pm e}=\xi_{-iu/2}^-\pm\tfrac12\Bracks1{\xi_e^-,\xi_{-iu/2}^-}$. Then 
\begin{align*}
	[\zeta_u^e,\zeta_v^{-e}](0)&=\tfrac18[[\xi_e^-,\xi_{iv}^-],\xi_{iu}^-](0)+\tfrac18[[\xi_e^-,\xi_{iu}^-],\xi_{iv}^-](0)\\
	&=\tfrac14\Parens1{\trip0e{(iv)}{(iu)}-\trip0{(iv)}e{(iu)}+\trip0e{(iu)}{(iv)}-\trip0{(iu)}e{(iv)}}\\
	&=\tfrac14\Parens0{2\cdot u\circ v+\trip0evu+\trip0euv}=u\circ v\ ,
\end{align*}
because $u^*=u$ and (\ref{eq:tripvscircprod}) give
\[
	\trip0euv=\trip0evu=e\circ(u^*\circ v)-u^*\circ(v\circ e)+v\circ(e\circ u^*)=u\circ v\ .
\]
\end{proof}

\begin{Lem}[nonredface-nilpotentpart]
	Under the assumptions of \thmref{Prop}{facestructure2}, $\ger g^e[\pm2]\cap F$ are faces of $F$, and $\ger g^e[-2,2]\cap F=\ger g^e[-2]\cap F\oplus\ger g^e[2]\cap F$. 
\end{Lem}

\begin{proof}
	Let $p_2^\pm$ and $p_2=p_2^++p_2^-$ be the orthogonal projections onto $\ger g^e[\pm2]$ and $\ger g^e[-2,2]$, respectively. By \thmref{Lem}{gradingcone}, $p^\pm_2(\ger z_F\cap F)\subset\ger g^e[\pm2]\cap F$, and the converse inclusion is obvious. By the same lemma, $\ger g^e[-2,2]\cap F=\ger g^e[-2]\cap F\oplus\ger g^e[2]\cap F$. 
	
	Now let $x,y\in F$ \scth $x+y\in\ger g^e[2]\cap F$. Write $x=\sum_{j=-2}^2x_j$, $y=\sum_{j=-2}^2y_j$ where $x_j,y_j\in\ger g^e[j]$. Then $x_{-2}+y_{-2}=0$, and with $x_{-2},y_{-2}\in F$, this implies $x_{-2}=y_{-2}=0$. By \thmref{Lem}{gradingcone} and \thmref{Prop}{nilpotent-cone}, $x_{-1}=y_{-1}=0$. Then $x_0+y_0=0$ where $x_0,y_0\in F$, and this implies $x_0=y_0=0$. Then $x_{-1}=y_{-1}=0$ (\emph{loc. cit.}). We conclude that $x=x_2,y=y_2\in\ger g^e[2]\cap F$, so this is a face. Analogously, $\ger g^e[-2]\cap F$ is a face.
\end{proof}

\begin{Lem}[tripotentpeirceinv]
	Let $c_1\perp c_2$ be tripotents. Then $Z_j(c_1+c_2)=Z_j(c_1-c_2)$ for $j=0,\tfrac12,1$. 
\end{Lem}

\begin{proof}
	Suffices to observe $\jBox{(c_1\pm c_2)}{(c_1\pm c_2)}=\jBox{c_1}{c_1}+\jBox{c_2}{c_2}$. 
\end{proof}

\begin{proof}[of \thmref{Prop}{facestructure2}]
	The spaces $\ger g^e[\pm1]$ are zero or faithful $\ger g_0(e)$-modules \cite[Chapter III, \S~4, Proposition~4.4, Corollary~4.5]{satake}. Hence, 
	\begin{equation}\label{eq:centraliser}
	\ger z_F=\ger g_F\cap(\ger m^e\oplus\ger g_1(e)\oplus\ger g^e[-2,2])
	\end{equation}
	where we recall $\ger g^e[0]=\ger g_0(e)\oplus\ger g_1(e)\oplus\ger m^e$ from (\ref{eq:zerograding-decomp}).

	Let $u^\pm\in X_1(e)$ \scth $\phi^{\pm e}(u^\pm)\in F$. Then $\phi^{\pm e}(u^\pm)\in\ger z_F$, and by \thmref{Lem}{commjordanprod}, $u^+\circ u^-=[\phi^e(u^+),\phi^{-e}(u^-)](0)=0$. On the other hand, $u^\pm\in\pm\Omega_1(e)$ by \thmref{Prop}{nilpotent-cone}, and $0=\Scp0e{u^+\circ u^-}=\Scp0{e\circ u^+}{u^-}=\Scp0{u^+}{u^-}$. If $c^\pm\sle e$ are \scth $\pm\Omega_1(c^\pm)$ is the face of $\Omega_1(e)$ generated by $u^\pm$, then $\Omega_1(c^+)\perp-\Omega_1(c^-)$. In particular, $c^+\perp c^-$. 

	Let $F^\pm=\ger g^e[\pm 2]\cap F$. By \thmref{Lem}{nonredface-nilpotentpart}, $F^\pm$ are nilpotent faces of $\Omega$, so by \thmref{Th}{mincone-co-orbit}, there exist tripotents $c^\pm\sle e$ \scth $F^\pm=\phi^{\pm e}(\Omega_1(\pm c^\pm))$. Necessarily, $c^+\perp c^-$. By \cite[Corollary~5.12]{loos_boundsymmdom}, there exists some $\ell\in K$ \scth $\ell(\pm c^\pm)=c^\pm$, and $\ell(e-(c^++c^-))=e-(c^++c^-)$. 

	By the above considerations, whenever $X_1(c)\subset X_F^\pm$, then $X_F^\mp\subset X_1(e-c)$. There exist $c^\pm\sle e$ \scth $X_1(c^\pm)\subset X_F^\pm$ (e.g.~$0=X_1(0)\subset X_F^+$). Then $c^+\perp c^-$, and $X_F^+\times X_F^-\subset X_1(e-c^-)\times X_1(e-c^+)$. By \thmref{Lem}{tripotentpeirceinv}, the Peirce decompositions for the tripotents $e$ and $\ell(e)$ are identical. On the other hand, it is clear by (\ref{eq:kaction-nilpotent}) that $\Ad(\ell)(\ger g^{c_-}[-2])=\ger g^{c_-}[2]\subset\ger g^e[2]$.

	Thus, if we set $F'=\Ad(\ell)(F)$, then we obtain $\ger g_{F'}=\ger g_0(\ell(e))\oplus\ger n_{F'}$ where the nilradical $\ger n_{F'}=[\ger n_{F'},\ger t_0(\ell(e))]\oplus\ger z_{F'}$, and 
	\[
		\ger z_{F'}\cap\ger g^{\ell(e)}[-2,2]=\Ad(\ell)(F\cap\ger g^e[-2,2])=\ger g^{c_+}[2]\oplus\ger g^{c_-}[2]\subset\ger g^e[2]\ .
	\]
	Furthermore, $\ger g^e[2]\cap F'$ is a nilpotent face by \thmref{Lem}{tripotentpeirceinv}, and therefore equals $\Omega_1(c)$ \fs $c\sle e$, by \thmref{Th}{mincone-co-orbit}. In particular, $\ger g^e[2]\cap F'=\Ad(\ell)(\ger g^e[-2,2]\cap F)$ is an irreducible cone. Hence, one of the faces $\ger g^e[\pm2]\cap F$ must be trivial. 

	Possibly replacing $e$ by $-e$ (which does not change $F$ or $\ger g_0(e)$), we may assume that $\ger g^e[-2]\cap F=0$. We set $c=c^+$. Arguing as usual with \thmref{Lem}{gradingcone} and \thmref{Prop}{nilpotent-cone}, we find that $\ger g_F\subset\ger q^e=\ger g^e[0,1,2]$, so that $\ger g_0(e)\oplus\ger g^c[2]\subset\ger g_F\subset\ger g_0(e)\ltimes\ger h^e$, by \thmref{Prop}{cone-gradingzero}. Moreover, $\ger z_F\subset\ger g_1(e)\oplus\ger m^e\oplus\ger g^c[2]$ by (\ref{eq:centraliser}). It follows that $\ger z_F=\ger g^c[2]$, $F\cap\ger z_F=\Omega_1(c)$, and $\ger n_F\subset\ger h^e$. 
\end{proof}

\begin{Prop}[nonredface]
	Let $\Omega=\Omega^\pm$ and $F\subset\Omega$ be a face. Assume that $\ger g_F=\ger g_0(e)\ltimes\ger n_F$ is not reductive. Then $\rk e<r$. Possibly replacing $e$ by $-e$, we have $\ger z_F=\ger g^c[2]=X_1(c)$ for a unique $c\sle e$, $\ger g_F=\ger g_0(e)\ltimes\ger h^{e,c}$ where 
	\[
		\ger h^{e,c}=\Set1{\eta_u^e}{u\in Z_{1/2}(e)\cap Z_{1/2}(c)}\oplus\ger g^c[2]\ .
	\]
	In particular, we have $F=F^\pm_{e,c}=F_e^\pm\cap F_c^\pm$, and this is an exposed face of $\Omega$. 
\end{Prop}

In addition to \thmref{Prop}{facestructure2}, the \emph{proof}Êrequires only the following simple lemma.

\begin{Lem}[heisenberg-subalg]
	Let $c\sle e$, and $\ger h\subset\ger h^e$ a subalgebra such that $\ger h\cap\ger g^e[2]=\ger g^c[2]$. Let $I$ be the complex structure on $\ger g^e[1]$ induced by that of $Z_{1/2}(e)$. If $\ger h\cap\ger g^e[1]$ is $I$-invariant, then $\eta_u^e\in\ger h$ implies $u\in Z_{1/2}(e)\cap Z_{1/2}(c)$. 
\end{Lem}

\begin{proof}
	Let $\eta_u^e\in\ger h$. Then $\ger g^c[2]\ni[I\eta_u^e,\eta_u^e]=\zeta_{-iv/2}^e$ where $v=q_e(iu,u)=8\trip0uue$ by (\ref{eq:heisenbergbracket}). By \thmref{Prop}{heisenberg-alg}, if $u\neq0$, then $v\neq0$. In particular, $v\in X_1(c)\setminus0$. 
	
	Now, $Z_{1/2}(e)=Z_{1/2}(e)\cap Z_{1/2}(c)\oplus Z_{1/2}(e)\cap Z_0(c)$. If $a$ lies in the second summand, then $\trip0aae=\trip0aa{(e-c)}\in Z_0(c)$. Similarly, $\trip0abe\in Z_{1/2}(c)$ if $a$ lies in the first summand, and $b$ lies in the second. Because $h_e(a,b)=8\trip0abe$ is $\Omega_1(e)$-positive Hermitian by \thmref{Prop}{heisenberg-alg}, $\trip0abe^*=\trip0bae$, and we conclude that $v\in X_1(c)$ if and only if $u\in Z_{1/2}(e)\cap Z_{1/2}(c)$. 
\end{proof}

\begin{proof}[of \thmref{Prop}{nonredface}]
	If we had $\ger g_0(e)=0$, then $\ger g_F$ would be nilpotent and hence Abelian \cite[Lemma~I.13]{hno_conalheisenberg}. By the assumption, this is excluded, so $\ger g_0(e)\neq0$. Let $\ger t$ be chosen according to (\ref{eq:cptcsa}) for some frame adapted to $e$, and $\ger t_F=\ger t_0(e)\oplus\ger z_F$ be the associated compact Cartan subalgebra of $\ger g_F$. Since $\ger g_F$ is not reductive, we have $\ger n_F=[\ger n_F,\ger t_0(e)]\oplus\ger z_F$ and the first summand contains no $\ger g_0(e)$-fixed vector \cite[Theorem~V.1]{neeb-classcones}. Hence, $\ger z_F=\ger z_{\ger g_F}(\ger g_0(e))$. 

By \thmref{Prop}{facestructure2}, possibly replacing $e$ by $-e$, we have $\ger g_F\subset\ger q^e$, $\ger n_F\subset\ger h^e$, and there exists a unique tripotent $c\sle e$ \scth $\ger z_F=\ger g^c[2]$ and $F\cap\ger z_F=\Omega_1(c)$. 

Hence, we have $\ger n_F\cap\ger g^e[2]=\ger g^c[2]$. On the other hand, $h_0'=iz\frac\partial{\partial z}-i\jBox ee\in\ger k_0(e)$ and $h_0'(z)=\frac i2z$ \fa $z\in Z_{1/2}(e)$. By (\ref{eq:g0g1action}), it follows that the $\ger g_0(e)$-module $\ger n_F\cap\ger g^e[1]$ is invariant under the complex structure $I$ of $\ger g^e[1]$. Invoking \thmref{Lem}{heisenberg-subalg}, it follows that $\eta_u^e\in\ger n_F$ implies $u\in Z_{1/2}(e)\cap Z_{1/2}(c)$. Since $\ger n_F=[\ger t_0(e),\ger n_F]\oplus\ger z_F$, we deduce that $\ger n_F=\ger n_F\cap\ger g^e[1]\oplus\ger g^c[2]\subset Z_{1/2}(e)\cap Z_{1/2}(c)\oplus\ger g^c[2]$. 

Let $\ger l=\ger g_0(e)\ltimes(Z_{1/2}(e)\cap Z_{1/2}(c)\oplus\ger g^c[2])$. Then $\ger h=\ger t_0(e)\oplus\ger g^c[2]$ is a compact Cartan subalgebra of $\ger l$. We have $\ger g_F\subset\ger l$ and $F_{e,c}^\pm=F_e^\pm\cap F_c^\pm\subset\ger l$ since $\ger g_0(e)\subset\ger g_0(c)$, $\ger g^e[2]\supset\ger g^c[2]$, $\ger g_0(c)\cap\ger g^c[2]=0$, and by the argument in the previous paragraph. The face $F_{e,c}^\pm$ is $\ger l$-invariant since it is $(\ger g_0(e)\ltimes\ger h^e)\cap(\ger g_0(c)\ltimes\ger h^c)$-invariant, and we have $\ger l\cap\ger h=\omega_0^\pm(e)\times\Omega_1(c)=F\cap\ger h$. It follows that $F\subset F_{e,c}^\pm$ and $F_{e,c}^\pm$ is regular in $\ger l$ (by the same argument as in the proof of \thmref{Prop}{cone-gradingzero}). 

But since $F_{e,c}^\pm\cap\ger h$ contains an element of the relative interior $F_{e,c}^{\pm\circ}$ \cite[Proposition III.5.14 and proof]{hilgert-hofmann-lawson}, it follows that the faces $F$ and $F_{e,c}^\pm$ are identical. In particular, $\ger g_F=\ger l$, and since the lattice of exposed faces is complete, $F$ is an exposed face. 
\end{proof}

\begin{Cor}[redfacealg]
	Let $F\subset\Omega=\Omega^\pm$ be a face with reductive face algebra $\ger g_F$. Then $F$ is a semi-simple face of the form $F=\Omega_0^\pm(e)$, or $\ger g_F$ is Abelian.
\end{Cor}

\begin{proof}
	By assumption and \thmref{Prop}{face-semisimplepart}, $\ger g_F=\ger g_0(e)\oplus\ger z_F$. We may assume that $\rk e<r$ since otherwise $\ger g_F=\ger z_F$ is Abelian. Then \thmref{Prop}{facestructure2} implies that (after possibly replacing $e$ by $-e$) there exists a tripotent $c\sle e$ \scth $\ger z_F=\ger g^c[2]$. We may assume $c>0$ since otherwise $F=\Omega_0^\pm(e)$. But then $F\supset\Omega_0^\pm(e)\oplus\Omega_1(c)$ and the latter cone contains points in the relative interior of $F_{e,c}^\pm$ \cite[Proposition III.5.14 and proof]{hilgert-hofmann-lawson}. Since $F\subset(\ger g_0(e)\oplus\ger g^c[2])\cap\Omega\subset F_{e,c}^\pm$ and is face, we conclude $F=F_{e,c}^\pm$. But this is a contradiction, since $F_{e,c}^\pm$ spans a non-reductive subalgebra of $\ger g$. 
\end{proof}

\subsection{Exhaustion of the faces of $\Omega^-$}

We are finally ready to describe all the convex faces of $\Omega^-$. 

\begin{Lem}[face-jordan]
	Let $F\subsetneq\Omega$ be a proper face. For any $\xi\in\ger g$, denote its Jordan decomposition by $\xi=\xi_s+\xi_n$. For all $\xi\in F^\circ$, we have $\xi_s,\xi_n\in F$. 
\end{Lem}

\begin{proof}
We have $\xi_s\in\Omega$ and $\xi_n\in\Omega^-\subset\Omega$ \fa $\xi\in\Omega$ \cite[Lemma IV.4]{neeb-convcoadj}. Let $\xi\in F^\circ$. Elements of $\Omega^\circ$ are elliptic and hence semi-simple, and $\partial\Omega$ is closed. Thus, $\xi_s+t\xi_n,t\xi_s+\xi_n\in\partial\Omega$ \fa $t\sge0$. This means that the line segments $[\xi,\xi_s]$ and $[\xi,\xi_n]$ lie within a proper face of $\Omega$, and therefore the open segments intersect $F^\circ$. But this implies $[\xi,\xi_s],[\xi,\xi_n]\subset F$. Hence, $\xi_s,\xi_n\in F$.
\end{proof}

\begin{Lem}[extremalray-invcone]
	Let $F$ be an extreme ray of $\Omega^\pm$. Then $F\subset\Omega^-$ if and only if $F$ is nilpotent. In this case, $F=\reals_{\sge0}\cdot X_c^+$ \fs primitive tripotent $c$. If this is not the case, then $F$ is conjugate to an extreme ray of $\omega^+$ which is not contained in $\omega^-$.
\end{Lem}

\begin{proof}
	Let $\xi\in F^\circ$, $\xi=\xi_s+\xi_n$. Then $\xi_s,\xi_n\in F=\reals_{\sge0}\cdot\xi$ by \thmref{Lem}{face-jordan}, and $\xi$ is semi-simple or nilpotent. The case of $\xi$ nilpotent is covered by \thmref{Th}{mincone-co-orbit}. Since $\pm X_e^\pm\in\Omega^-$ and $i\cdot\jBox ee=\tfrac12(X_e^+-X_e^-)$ for any tripotent $e$, no extreme ray of $\omega^-$ is an extreme ray of $\Omega^\pm$, by \thmref{Cor}{extremalgens-tripotents}. 
	
	Hence, $\xi$ is semi-simple if and only if $F$ is an extreme ray of $\Omega^+$ which is not contained in $\Omega^-$. In this case, $\xi$ is conjugate to an element of $\omega^+$ by \thmref{Prop}{cone-regular}, so we may assume $F\subset\omega^+$. Since $\omega^+\subset\Omega^+$, $F$ is then an extreme ray of $\omega^+$. 
\end{proof}

\begin{Cor}[solvfacealg]
	Any face of $\Omega^-$ with a solvable face algebra is a nilpotent face. 
\end{Cor}

\begin{proof}
	Let $F\subset\Omega^-$ be a face with solvable face algebra. By \thmref{Prop}{face-structure}, $\ger g_F$ is Abelian. By Strasziewicz's spanning theorem, the cone spanned by the extreme rays of $F$ is dense in $F$. Hence, there exists $x\in F^\circ$ which is the positive linear combination of extreme generators. By \thmref{Lem}{extremalray-invcone}, all of the latter are nilpotent elements of $\ger g$. Since the commute, $x$ is also nilpotent, and $F$ is by definition a nilpotent face. 
\end{proof}

\begin{Prop}[semsimface]
	Let $F$ be a semisimple face of $\Omega=\Omega^\pm$. Then $F=\Omega_0^\pm(e)$ \fs tripotent $e$, or $\ger g_F$ is Abelian and contained in a compact Cartan subalgebra of $\ger g$. The latter alternative only occurs for $F=0$ or $\Omega=\Omega^+$, and in this case, $F$ is conjugate to a face of $\omega^+$. In particular, the set of semi-simple faces of $\Omega$ is a lattice. 
\end{Prop}

\begin{proof}
	By the semisimplicity of $F$, $F^\circ\subset\ger g\setminus0$ contains semi-simple elements. By \thmref{Th}{mincone-co-orbit}, the nilpotent faces consist of nilpotent elements of $\ger g$. Hence, $F$ is not contained in a nilpotent face. But then $F^\circ$ cannot intersect any nilpotent face. Since any nilpotent element of $\Omega^+$ is contained in a nilpotent face, $F^\circ$ consists of semi-simple elements. 
	
	Hence, we may choose $\ger t$ according to (\ref{eq:cptcsa}) \scth $\ger t_F=\ger t\cap\ger g_F$ is a compact Cartan subalgebra of $\ger g_F$. In particular, $\ger g_F$ is $\ger t$-invariant. There exists an additively closed subset $P\subset\Delta$ \scth $\ger g_{F\cplxs}=\ger g_{\ger t_F,P}$ where $\ger g_{\ger h,P}=\ger h_\cplxs\oplus\bigoplus_{\alpha\in P}\ger g_\cplxs^\alpha$ \fa $\ger h\subset\ger t$. Let $\ger t'_F=\bigcap_{\alpha\in P\cap(-P)}\ger t_F\cap\ker\alpha$ and $Q=P\setminus(-P)$. Then the (nil) radical of $\ger g_{F\cplxs}$ is $\ger r_{F\cplxs}=\ger n_{F\cplxs}=\ger g_{\ger t'_F,Q}$, and $[\ger n_{F\cplxs},\ger n_{F\cplxs}]=\ger g_{0,S}$ where
	\[
		S=(Q+Q)\cap\Delta\cup\Set0{\alpha\in Q}{\alpha(\ger t'_F)\neq0}\ .
	\]
	All of this follows from \cite[Chapter VIII, \S~3.1, Proposition 2]{bourbaki}. On the other hand, $[\ger n_F,\ger n_F]$ must be central in $\ger g_F$, by \thmref{Prop}{face-structure}. In particular, $[\ger n_F,\ger n_F]\subset\ger t_F$. This implies $S=\vvoid$ and $[\ger n_F,\ger n_F]=0$. Then $\ger n_F$ is an Abelian ideal of $\ger g_F$, and therefore central \cite[Proposition~VII.3.23]{neeb_holconv}. It follows that $Q=\vvoid$, and $\ger n_F=\ger t'_F$. In particular, $\ger g_F$ is reductive.
	
	It follows from \thmref{Prop}{face-semisimplepart} that $\ger g_F=\ger g_0(e)\oplus\ger t_F'$ \fs tripotent $e$, where $\ger t_F=\ger t_0(e)\oplus\ger t_F'$. Since $\ger t_F\cap\Omega^\pm=\omega_0^\pm(e)\oplus\ger t_F'\cap\omega^\pm$, we have $F=\Omega_0^\pm(e)\oplus\ger t_F'\cap\omega^\pm$. By \thmref{Cor}{redfacealg}, $F=\Omega_0^\pm(e)$ or $F=\ger t_F'\cap\omega^\pm$. The latter alternative is impossible for $\Omega=\Omega^-$ in view of \thmref{Cor}{solvfacealg}, since $\ger t_F'$ consists of semi-simple elements. 
\end{proof}

We summarise our considerations for $\Omega^-$ in the following theorem and corollary.

\begin{Th}[facemincone]
	Each face of $\Omega^-$ is one of $\Omega_0^-(e)$, $\Omega_1(e)$, $e$ a tripotent; or of $F_{e,c}^-$, $e\sge c>0$ tripotents with $\rk e<r$. In particular, $\Omega^-$ is facially exposed. 
\end{Th}

\begin{proof}
 	First we remark that $F_{e,0}^-=\Omega_0^-(e)$ for $\rk e<r$ and that $F_{e,c}^-=\Omega_1(c)$ for $\rk e=r$. Thus, if $F$ is a face with $\ger g_F$ non-reductive, then $F$ occurs in the second part of the list set out above, by \thmref{Prop}{nonredface}. If $\ger g_F$ is reductive, then by \thmref{Cor}{redfacealg}, \thmref{Prop}{semsimface}, \thmref{Cor}{solvfacealg} and \thmref{Th}{mincone-co-orbit}, $F$ occurs in the first part of the list. 
\end{proof}

\begin{Cor}[facemincone]
	Every face of $\Omega^-$ is one of the faces $F_{e,c}^-$, for tripotents $0\sle c\sle e$. 
\end{Cor}

\paragraph{Conjugacy classes of faces and $K$-orbit type decomposition of $\Omega^-$}

\begin{Th}[faceconj]
	Any two faces $F_{e,c}^\pm$ and $F_{e',c'}^\pm$ of $\Omega^\pm$ are $G$-conjugate if and only one has $(\rk e,\rk c)=(\rk e',\rk c')$, if and only if they are $K$-conjugate. In particular, 
	\[
		\Omega_{k,\ell}=\bigcup\nolimits_{c\sle e\,,\,(\rk e,\rk c)=(k,\ell)}F_{e,c}^{-\circ}\qquad,\qquad0\sle\ell\sle k\sle r
	\]
	are exactly the orbit types of the $K$-action on $\Omega^-$. If $M_{k,\ell}$ is the set of pairs $(e,c)$ of tripotents $e\sge c$ \scth $(\rk e,\rk c)=(k,\ell)$, then $\Omega_{k,\ell}$ is a $K$-equivariant fibre bundle over the $K$-homogeneous space $M_{k,\ell}$ with typical fibre $F_{e,c}^{-\circ}$. 
\end{Th}

\begin{proof}
	Given the equality of ranks, the faces are $K$-conjugate, in view of \cite[Theorem 5.9]{loos_boundsymmdom}. Moreover, they are certainly $G$-conjugate if they are $K$-conjugate. If they are $G$-conjugate, then the algebras $\ger g^c[2]$ and $\ger g^{c'}[2]$ are $G$-conjugate, and so are $\ger g_0(e)$ and $\ger g_0(e')$, as the centres of the respective face algebras, and their Levi complements invariant under compact Cartan subalgebras, respectively. By \thmref{Th}{mincone-co-orbit}, we have $\rk c=\rk c'$, and by \thmref{Lem}{g0conj}, we have $\rk e=\rk e'$. 
	
	Any element in the relative interior of the face $F=F_{e,c}^-$ is $G$-conjugate to an element of the relative interior of $f=F\cap(\ger t_0(e)\oplus\ger g^c[2])=\omega_0^-(e)\oplus\Omega_1(c)$. Moreover, $k\in K$ fixes $\xi\in f^\circ$ if and only if $k$ fixes $\xi_s\in\omega_0^-(e)^\circ$ and $\xi_n\in\Omega_1(c)^\circ$ where we denote by $\xi=\xi_s+\xi_n$ the Jordan decomposition. The stabiliser of $\xi_s$ is $N_K(\ger t_0(e))$, and the stabiliser of $\xi_n$ is $K^c$, independent of $\xi$. This shows that $\Omega_{k,\ell}$, $(k,\ell)=(\rk e,\rk c)$, is exactly a single $K$-orbit type. By \thmref{Cor}{facemincone}, the assertion follows. 
\end{proof}

\begin{Cor}
	For any $r\sge k\sge\ell\sge 0$, $\Omega_{k,\ell}$ is $K$-equivariantly homotopy equivalent to the $K$-homo\-gen\-eous space $M_{k,\ell}=K/(K^e\cap K^c)$ (where $(e,c)\in M_{k,\ell}$). 
\end{Cor}

\section{The stratification of the minimal Ol'shanski\u\i{} semigroup}

In this section, we apply our previous considerations to achieve our ultimate goal: The decomposition of the minimal Ol'shanski\u\i{} semigroup into $K$-orbit type strata, and their description in terms of $K$-equivariant fibre bundles. 

\subsection{The minimal Ol'shanski\u\i{} semigroup}

	There exists a connected complex Lie group $G_\cplxs$ with Lie algebra $\ger g_\cplxs$ \scth $G\subset G_\cplxs$. (\emph{E.g.}, consider the projective completion $D^*$ of $D$ \cite[\S\S\,8--9]{loos_boundsymmdom} and let $G_\cplxs$ be the connected component of $\Aut0{D^*}$ \cite[Proposition 9.4]{loos_boundsymmdom}. Alternatively, we may invoke \cite[Proposition~25.9]{gloeckner_complex}.)
	
	By the following proposition, there exists a closed complex semigroup $\Gamma\subset G_\cplxs$ \scth $\Gamma=G\cdot\exp i\Omega^-$ and $G\times\Omega^-\to\Gamma:(g,\xi)\mapsto g\exp i\xi$ is a homeomorphism which restricts to a diffeomorphism $G\times\Omega^{-\circ}\to\Gamma^\circ$. This semigroup is called the \emph{minimal Ol'shanski\u\i{} semigroup}. 
	
	The following proposition is a compilation of known results. We give it for the reader's convenience, since we lack a succinct reference. The construction of Ol'shan\-ski\u\i{} semigroups is developed in full generality in \cite[Chapters~3, 7]{hn_liesemigrp}, \cite[Chapter~XI]{neeb_holconv}.

\begin{Prop}[olshanski]
Let $H\subset H_\cplxs$ be connected Lie groups where $H$ is closed, $H_\cplxs$ is complex, and the Lie algebra of $H_\cplxs$ is $\ger h_\cplxs$. Let $\Omega\subset\ger h$ be an invariant regular convex cone. Then $\psi:H\times\Omega\to H_\cplxs:(h,\xi)\mapsto h\cdot\exp i\xi$ is a homeomorphism onto a closed subsemigroup $\Gamma$, and induces a diffeomorphism $H\times\Omega^\circ\to\Gamma^\circ$ where $\Gamma^\circ\subset\Gamma$ is dense.
\end{Prop}

\begin{proof}
Let $\phi:\widetilde H_\cplxs\to H_\cplxs$ be the universal covering. Then $\ker\phi$ is a discrete central subgroup. The $\mathrm{Gal}(\cplxs:\reals)$-action on $H_\cplxs$ lifts to $\widetilde H_\cplxs$. Let $\widetilde H$ be the fixed group of this action; then $\widetilde H$ is closed and connected \cite[Chapter IV, Theorem 3.4]{loos_symmetric1}. The adjoint action of $\ger h$ has imaginary spectrum \cite[Proposition~II.2]{neeb_invconv}. 

Therefore, $\tilde\psi:\widetilde H\times\Omega\to\widetilde H^\cplxs:(h,\xi)\mapsto h\exp i\xi$ (where we take the exponential map of $\widetilde H_\cplxs$) is a homeomorphism onto a closed subsemigroup $\widetilde\Gamma\subset\widetilde H_\cplxs$ which restricts to a diffeomorphism $\widetilde H\times\Omega^\circ\to\widetilde\Gamma^\circ$ \cite[Theorems~XI.1.7, XI.1.10]{neeb_holconv}. In particular, $\widetilde H$ is a retraction of $\widetilde H_\cplxs$ and therefore simply connected. It follows that $H$ and $H_\cplxs$ are homotopy equivalent \cite[Proposition~25.9]{gloeckner_complex}. We have canonical isomorphisms
\[
\ker\phi\cap\widetilde H\to\pi_1(H,1)\to\pi_1(H_\cplxs,1)\to\ker\phi\ .
\]
This map associates to $h\in\ker\phi$ the homotopy class of $\phi\circ\gamma_h$, $\gamma_h$ a path in $\widetilde H$ from $1$ to $h$; to this, the homotopy class in $H_\cplxs$ of $\phi\circ\gamma_h$; hereto, the end point of a lifting of $\phi\circ\gamma_h$ in $H_\cplxs$. Since $\gamma_h$ is such a lifting and $\gamma_h(1)=h$, the composite map is the identity, and $\ker\phi\subset\widetilde H$. Thus, we conclude that $\tilde\psi$ drops to a map $\psi$ with the required properties. 
(For the statement on the interiors, see \cite[Lemma XI.I.9]{neeb_holconv}.) 
\end{proof}

\subsection{The stratification of the minimal Ol'shanski\u\i{} semigroup into $K$-orbit types}

\begin{Def}
	A \emph{Lie semigroup} is a pair $(S,H)$ where $H$ is a connected Lie group and $S\subset H$ is a closed subsemigroup which is generated (as a closed subsemigroup) by the one-parameter semigroups it contains \cite[Definition IV.3]{neeb-invsubsemgrps}. The \emph{tangent wedge} of $(S,H)$ is the convex cone $L(S)=\Set1{\xi\in\ger h}{\exp(\reals\xi)\subset S}$ \cite[Definition IV.2]{neeb-invsubsemgrps}.  

	Let $(S,H)$ be a Lie semigroup. A \emph{face} of $(S,H)$ is a subsemigroup $F\subset S$ \scth $S\setminus F$ is a semigroup ideal. 
\end{Def}

\begin{Prop}[olshanski-faces]
	The pair $(\Gamma,G_\cplxs)$ is a Lie semigroup whose only faces are $G$ and $\Gamma$. 
\end{Prop}

\begin{proof}
	Consider the cone $W=\ger g\oplus i\Omega^-$. It is $G$-invariant and therefore a \emph{Lie wedge} \cite[Definition IV.1]{neeb-invsubsemgrps}. It equals the tangent wedge of $\Gamma$ and is therefore \emph{global} \cite[Definition IV.23, Lemma IV.24]{neeb-invsubsemgrps}; in particular, $(\Gamma, G_\cplxs)$ is a Lie semigroup. By \cite[Lemma~7.30]{hn_liesemigrp}, \cite[Lemma~II.2.11]{hilgert-hofmann-lawson}, $W$ is \emph{Lie semialgebra} \cite[Definition~IV.29]{neeb-invsubsemgrps}. By \cite[Proposition~IV.32]{neeb-invsubsemgrps}, the faces of $(\Gamma,G_\cplxs)$ are among the closed subsemigroups whose tangent wedges are faces of $W$ and therefore of the form $\ger g\oplus iF$ where $F\subset\Omega^-$ is a faces. 
	
	Let $S\subset\Gamma$ be a non-trivial face. Then $L(S)=\ger g\oplus iF$ where $F\subset\Omega^-$ is a non-trivial face. Hence, $F$ contains an extreme ray: by \thmref{Lem}{extremalray-invcone}, $F$ intersects the minimal nilpotent orbit of $\Omega^-$ non-trivially. Since $G\subset S$, $L(S)$ is $G$-invariant, and therefore contains the minimal nilpotent orbit in $i\Omega^-$. Since $L(S)$ is a closed convex cone, we have $i\Omega^-\subset L(S)$ by \thmref{Cor}{minorbit-span}, and thus $L(S)=W$. This implies $\Gamma=S$. 
\end{proof}

	The stratification of $\Gamma$ into $K$-orbit types is more interesting. To describe it, let $F\subset\Omega^-$ be a face. Then $\ger g_F=\Span0F_\reals$ is a subalgebra, and we may consider the analytic subgroups $G_F\subset G$ and $G_{F\cplxs}\subset G_\cplxs$ associated with $\ger g_F$ and $\ger g_{F\cplxs}$, respectively. We have an Ol'shanski\u\i{} semigroup $\Gamma_F=G\cdot\exp iF\subset G_{F\cplxs}$ whose interior $\Gamma_F^\circ$ in $G_{F\cplxs}$ is $G\cdot\exp iF^\circ$ ($F^\circ$ denoting relative interior). 

\begin{Th}[olshanski-strata]
	The subsemigroups $\Gamma_F,\Gamma_{F'}\subset\Gamma$, $F=F_{e,c}^-$ and $F'=F_{e',c'}^-$, are $G$-conju\-gate if and only if they are $K$-conjugate, if and only if $(\rk e,\rk c)=(\rk e',\rk c')$. The orbit type strata for the action of $K$ on $\Gamma$ by conjugation are 
	\[
	\Gamma_{k,\ell}=G\cdot\exp i\Omega_{k,\ell}=\bigcup\nolimits_{(e,c)\in M_{k,\ell}}G\cdot\exp iF_{e,c}^{-\circ}\qquad,\qquad0\sle\ell\sle k\sle r\ .
	\]
\end{Th}

\begin{proof}
	This follows immediately from \thmref{Th}{faceconj}. 
\end{proof}

\begin{Cor}
	The orbit type stratum $\Gamma_{k,\ell}$ is $K$-equivariantly homotopy equivalent to $K\times M_{k,\ell}=K\times(K/(K^e\cap K^c))$ (where $(e,c)\in M_{k,\ell}$). 
\end{Cor}

\begin{proof}
	Clearly, $\Gamma_{k,\ell}$ fibres over $M_{k,\ell}$ with fibre $G\cdot\exp iF_{e,c}^{-\circ}$. Moreover, there exists a $K$-equivariant homotopy equivalence s $G\cdot\exp iF_{e,c}^{-\circ}\simeq G\simeq K$. 
\end{proof}

\emph{Acknowledgements}. This paper is a completely rewritten and expanded version of a part of my doctoral thesis \cite{alldridge-phd} under the supervision of Harald Upmeier. I wish to thank him for his support and guidance.

\end{document}